\documentclass[11pt]{amsart}
\usepackage{palatino}
\usepackage[all]{xy}
\usepackage{amssymb}
\usepackage{latexsym}
\usepackage{amscd}
\usepackage{color}
\input amssym.def
\input amssym

\newtheorem{theorem}{Theorem}[section]
\newtheorem{lemma}[theorem]{Lemma}
\newtheorem{definition}[theorem]{Definition}

\newtheorem{proposition}[theorem]{Proposition}
\newtheorem{corollary}[theorem]{Corollary}

\newtheorem{remark}{Remark}

\newcommand{\ds}{\displaystyle}

\newcommand{\halmos}{\rule{1ex}{1.4ex}}

\newcommand{\bea}{\begin{eqnarray}}
\newcommand{\epfv}{\hspace*{\fill}\mbox{$\halmos$}}
\newcommand{\nb}{{\mbox {\tiny{ ${\bullet \atop \bullet}$}}}}
\newcommand{\eea}{\end{eqnarray}}
\newcommand{\no}{\noindent}
\newcommand{\nn}{\nonumber \\}
\newcommand{\be}{\begin{equation}}

\newcommand{\ee}{\end{equation}}

\newcommand{\lo}{{\rm log}}

\usepackage{amssymb}

\begin{document}
\title[Logarithmic operators and vertex operators]{Logarithmic intertwining operators and vertex operators}
\author{Antun Milas}

\begin{abstract}
\noindent This is the first in a series of papers where we study
logarithmic intertwining operators for various vertex subalgebras
of Heisenberg and lattice vertex  algebras. In this paper we
examine logarithmic intertwining operators associated with rank
one Heisenberg vertex operator algebra $M(1)_a$, of central charge
$1-12a^2$. We classify these operators in terms of { depth} and
provide explicit constructions in all cases. Our intertwining
operators resemble { puncture operators} appearing in quantum
Liouville field theory. Furthermore, for $a=0$ we focus on the
vertex operator subalgebra $L(1,0)$ of $M(1)_0$ and obtain
logarithmic intertwining operators among indecomposable Virasoro
algebra modules. In particular, we construct explicitly a family
of { hidden} logarithmic intertwining operators, i.e., those that
operate among two ordinary and one genuine logarithmic
$L(1,0)$-module.

%Our work presents the first purely algebraic treatment of what
%physicists have been referring as logarithmic conformal field
%theory of central charge $c_p=1-\frac{6(p-1)^2}{p}$, $p \geq 1$.
%In the sequel \cite{AM}, \cite{M3} we will construct logarithmic
%intertwiners for various vertex subalgebras of $M(1)_a$, such as
%$\mathcal{W}(2,2p-1)$-algebras of central charge $c_p$, $p \geq
%2$.
\end{abstract}

\address{Department of Mathematics and Statistics,
University at Albany (SUNY), Albany, NY 12222}
\email{amilas@math.albany.edu}

\maketitle

\renewcommand{\theequation}{\thesection.\arabic{equation}}
\renewcommand{\thetheorem}{\thesection.\arabic{theorem}}
\setcounter{equation}{0} \setcounter{theorem}{0}
\setcounter{equation}{0} \setcounter{theorem}{0}
\setcounter{section}{-1}

\section{Introduction}

The theory of vertex algebras continues to be very effective in
proving rigorous results in two-dimensional Conformal Field Theory
(CFT) (for some recent breakthrough see \cite{H}, \cite{H1},
\cite{Le}).

In 1993, Gurarie \cite{Gu} studied a CFT-like structure with two
features absent in the ordinary CFT: logarithmic behavior of
matrix coefficients and appearance of indecomposable
representations of the Virasoro algebra underlying the theory.
%Such a structure contradicts to some very definitions of conformal
%field theory.
There are several additional examples of "logarithmic" models that
have been discovered since then (see for instance \cite{GK},
\cite{F1}, \cite{G1}, and especially \cite{FFHST}, \cite{F2},
\cite{G2} and references therein). By now, a structure that
involves a family of modules for the chiral algebra, closed under
the "fusion", with logarithmic terms in the operator product
expansion is usually called a {\em logarithmic conformal field
theory} (LCFT). The most important examples of LCFTs are the
so-called {\em rational} LCFTs. These involve only finitely many
inequivalent irreducible representations, but also some
indecomposable logarithmic representations so that the modular
invariance is preserved
%formulation of rational
%LCFT. Another important problems is to formulate an analogue of
%the Verlinde formula (if there is such).
%There are
%several proposals in the literature which explain how the
%logarithmic behavior of correlation functions is linked to
%nontrivial Jordan blocks under the action of $L(0)$ (cf.
%\cite{G2}, \cite{F2}).
(e.g., the triplet model \cite{GK}, \cite{GK1}, \cite{GK2},
\cite{BF}, etc.). We stress here that neither LCFT nor rational
LCFT are mathematically precise notions.

In \cite{M1} we proposed a purely algebraic approach to LCFT based
on the notion of logarithmic modules and {\em logarithmic
intertwining operators}. The key idea is to introduce a
deformation parameter $\lo(x)$ and to define logarithmic
intertwining operators as expressions involving intertwining-like
operators multiplied with appropriate powers of $\lo(x)$, such
that the translation invariance is preserved. These operators can
be used to explain appearance of logarithms in correlation
functions. In our setup we do not require an extension of the
space of "states" for the underlying vertex algebra. There are
other proposals in the literature \cite{FFHST} where $\lo(x)$ is
also viewed as a deformation parameter, but with an important
difference that $\lo(x)$ is also part of an extended chiral
algebra (or an OPE algebra). Even though the construction in
\cite{FFHST} has been shown successful in explaining various
logarithmic behaviors of CFTs, it is unclear to us if the approach
in \cite{FFHST} can be used to address the problem of fusion.
%Another approach similar in
%spirit, is proposed in \cite{FSST}.
%Furthermore and explains
%features observed by physicists on several examples.
%The idea is
%to threat $\lo(x)$ as a special deformation parameter (so $\lo(x)$
%is a formal variable), but with an important
%Our approach is different from some formulations in the literature
%where $\lo(x)$ is viewed as a deformation parameter of the chiral
%algebra. We keep the definition of the vertex operator algebra
%unchanged and instead deform the intertwining operator map.
On the other hand, logarithmic intertwining operators have been
used by Huang Lepowsky and Zhang in \cite{HLZ} as a convenient
tool to develop a generalization of Huang-Lepowsky's tensor
product theory \cite{HL} to non-semisimple tensor categories.
Another important contribution is Miyamoto's  generalization of
Zhu's modular invariance theorem for vertex operator algebras
satisfying the $C_2$-condition \cite{My}, which possibly involve
logarithmic modules. In view of \cite{My} we tend to believe that
rational LCFT give rise from vertex algebras satisfying the
$C_2$-cofiniteness condition.
%all rational LCFT arise from (not necessarily rational) VOAs
%satisfying the $C_2$-cofiniteness condition.
%But still, by now there is no consistent
%%
%We should say here
%that some portions of the logarithmic constructions in \cite{HLZ}
%and \cite{My} were discovered independently by the author.

From everything being said it appears that several important
aspects of LCFTs can be studied in the framework of vertex
(operator) algebras.

%Rational LCFT should correspond to VOA which satisfy
%$C_2$-cofiniteness property.
This paper continues naturally on \cite{M1} and \cite{M2}. In the
present work we focus on a simple, yet interesting class of vertex
operator algebras-those associated with Heisenberg Lie algebras
and the Virasoro algebra. The aim here is to construct a family of
logarithmic intertwining operators associated with certain weak
$M(1)_a$-modules, which can be used for building logarithmic
intertwining operators among indecomposable representations of the
Virasoro algebra and various $\mathcal{W}$-algebras ( \cite{AM},
\cite{M3}). The most interesting part of our work is an explicit
construction of the so-called hidden logarithmic intertwining
operators, i.e., those which intertwine a pair of ordinary and one
logarithmic module. We should stress here that our logarithmic
intertwining operators are closely related to {\em puncture
operators} appearing in quantum Liouville field theory \cite{Se},
\cite{ZZ}.

Let us elaborate the construction on an example. Consider the
Feigin-Fuchs module $M(1,\lambda)_a$ of central charge $c=1-12a^2$
and lowest conformal weight $\frac{\lambda^2}{2}-a \lambda$.
Tensor the module $M(1,\lambda)_a$ with a two-dimensional space
$\Omega$, where $h(0)$ acts on $\Omega$ (in some basis) as
$$\left[ \begin{array}{cc} \lambda & 1 \\ 0 & \lambda \end{array}
\right]$$ so we obtain a weak $M(1)_a$-module $M(1,\lambda)_a
\otimes \Omega$. Now it is easy to see that $M(1,\lambda)_a
\otimes \Omega$ is an ordinary module if and only if $\lambda=a$,
so that $M(1,a)_a \otimes \Omega$ is of lowest conformal weight
$-\frac{a^2}{2}$. Our main results (cf. Theorem \ref{uppertheorem}
and Theorem \ref{pre-upper}) then provides us with a genuine
logarithmic intertwining operator of type
$${ W \choose M(1,a)_a \otimes \Omega \ \ M(1,a)_a \otimes
\Omega},$$ where $W$ is a logarithmic module. Notice that for
special $a$ this result is in "agreement" with some results in the
physics literature. For instance, for $c=-2$ and $a=\frac{1}{2}$,
the lowest conformal weight of $M(1,\frac{1}{2})_{\frac{1}{2}}
\otimes \Omega$  is $\frac{-1}{8}$. This model is known to be
logarithmic after \cite{Gu}. More generally, one takes
$c=1-\frac{6(p-1)^2}{p}$, $p \geq 2$ (cf. \cite{AM}, \cite{M3}).

We should say here that Feigin-Fuchs modules $M(1,a)_a$ with the
lowest conformal weight $-\frac{a^2}{2}$ and central charge
$1-12a^2$ are indeed special from at least two points of view.
These modules are self-dual (cf. (\ref{cont-iso})) and in addition
we have
$$h-\frac{c}{24}=\frac{-1}{24},$$ so that
$${\rm tr}|_{M(1,a)_a} q^{L(0)-c/24}=\frac{1}{\eta(\tau)},$$
where $\eta(\tau)=q^{1/24} \prod_{i=1}^\infty (1-q^i)$ is the
Dedekind $\eta$-function. In fact these are the only Feigin-Fuchs
modules whose modified graded dimensions are modular functions.

%We should mention here that the intertwining operators we
%constructed are related to {\em puncture operators} appearing in
%quantum Liouville field theory \cite{MS}, \cite{ZZL}.
%genuine logarithmic module and we can con$c=1$ we obtain $a=0$,
%Our construction can be easily modified so that it applies to
%other $W$-algebras.
%approaches and interpretation of logarithmic conformal field
%theory LCFT (see \cite{M1}  and \cite{HLL}) one algebraic was
%addressed in \cite{M1}. It appears that none of the approaches can
%be used But two features persist which are intimately related:
%appearance of Jordan blocks under the action of the Virasoro
%generator $L(0)$ and the correlation functions have logarithmic
%singularities. In fact the two are closely \cite{M1}. One
%algebraic approach was proposed in \cite{M1}. Even though the
%proposal in \cite{M1} does not cover all known LCFT it is good
%enough to construct tensor categories \cite{HLL}. In our modest
%other approaches the only consistent approach is the one in
%\cite{M1}, where logarithms are constructed as deformations of the
%ordinary modules.

Let us briefly outline the content of the paper. In Section 2 we
recall the notions of logarithmic module and of logarithmic
intertwining operator. In sections 3 and 4 we prove some standard
results about extended Heisenberg Lie algebras, the corresponding
vertex operator algebras and associated logarithmic intertwining
operators. In Section 5 we derive a logarithmic version of the
Li-Tsuchiya-Kanie's "Nuclear Democracy Theorem" for logarithmic
intertwining operators needed for our construction. In sections 5
and 7 we provide a classification of logarithmic intertwining
operators among a triple of logarithmic $M(1)_a$-modules. Finally,
in sections 9 and 10 we use certain restriction theorems for
construction of logarithmic intertwining operators among triples
of logarithmic $L(1,0)$-modules. In particular, we obtain a family
of hidden intertwining operators from a pair of ordinary
$L(1,0)$-modules. Section 11 is an introduction to \cite{AM} and
\cite{M3}.

\vskip 2mm

This paper is based on a talk the author gave at the AMS Special
Session on Representations of Infinite Dimensional Lie Algebras
and Mathematical Physics in Bloomington, IN, April 2003. Because
of the recent increasing interest in LCFT we decided to make our
paper public.

\vskip 2mm

\noindent {\bf Acknowledgments:} We are grateful to anonymous
referees for several useful comments.

\section{Logarithmic modules and logarithmic intertwining operators}

In this section notation and definitions are mostly from
\cite{DLM}, \cite{HLZ}, \cite{LL} and \cite{M1}. For an accessible
introduction to vertex (operator) algebras and their
representations we refer the reader to \cite{LL}.

\begin{definition} \label{def-logmod}
{\em Let $(V,Y,{\bf 1},\omega)$, $V=\ds{\coprod_{i \in
\mathbb{Z}}} V_i$, be a vertex operator algebra. A {\em
logarithmic} $V$-module $M$ is a weak $V$-module \cite{DLM} which
admits a decomposition
$$M=\coprod_{ h \in \{ h_i +\mathbb{Z}_+, \  i=1,...,k\}} M_h,  \ \ h_i \in \mathbb{C}, \ \ {\rm for \ some} \ k \in \mathbb{N}$$
$$M_{h}=\{ v \in M : (L(0)-hI)^m v=0, \ {\rm for} \ \ {\rm some}
\ \ m \}, \ \ {\rm dim}(M_h) < +\infty.$$
%such that for every $v
%\in V_i$ and $m \in \mathbb{Z}$ we have
%$$ v_m M_{h} \subseteq M_{h+i-m-1}.$$
Here $M_h$ denotes the {\em generalized homogeneous subspace} of
$M$ of {\em generalized weight} $h$.
%such that for every $v \in M_h$ there exists $k \in \mathbb{N}$
%such that $(L(0)-hI)^k v=0$.
}
\end{definition}
It is important to mention that we use the definition of weak
module as in \cite{DLM}, so that $M$ carries an action of the full
Virasoro algebra and not only the $L(-1)$-operator. Hence, for
every $v \in V_i$ and $m \in \mathbb{Z}$ we have
$$ v_m M_{h} \subseteq M_{h+i-m-1}.$$
The category of logarithmic $V$-modules $\mathcal{LOG}$ is
essentially just a different name for the category $\mathcal{O}$
introduced in \cite{DLM}. In particular every logarithmic module
is {\em admissible}  \cite{DLM}. We say that a weak $V$-module is
a {\em generalized logarithmic module} if it decomposes into (not
necessarily finite-dimensional) generalized $L(0)$-eigenspaces.
Furthermore, we will say that a logarithmic $V$-module $M$ is {\em
genuine} if under the action of $L(0)$ it admits at least one
Jordan block of size $2$ or more.
%If the size of every Jordan block in is bounded
%by above we say that $M$ is {\em strongly logarithmic}.

We recall the definition of logarithmic intertwining operators
from \cite{M1}. More precisely, here we are using a slightly more
general definition following \cite{HLZ}. In what follows $\lo(x)$
is just a formal variable satisfying $\frac{d}{dx}
\lo(x)=\frac{1}{x}$.

\begin{definition} \label{def-logint}
{\em A logarithmic intertwining operator among a triple of
logarithmic $V$-modules $W_1$, $W_2$ and $W_3$ is a linear map
\bea \label{map} && \mathcal{Y}( \cdot, x) \cdot  : W_1 \otimes
W_2 \longrightarrow W_3 \{x\} [\lo(x)] \nn && \mathcal{Y}(w_1,
x)w_2=\sum_{n \in \mathbb{N}} \sum_{\alpha \in \mathbb{C}}
\lo^n(x) x^{-\alpha-1} (w_1)^{(n)}_\alpha w_2, \ \
(w_1)^{(n)}_\alpha \in {\rm End}(W_2,W_3), \eea satisfying:
\begin{itemize}
\item[(i)] (Truncation condition) For every $w_1 \in W_1$, $w_2
\in W_2$ and $\alpha \in \mathbb{C}$,
$$(w_1)_{\alpha+m}^{(k)} w_2=0, \ \ {\rm for} \ {m \in
\mathbb{N}}, \ m >>0.$$
 \item[(ii)]
(Translation invariance) For every $w_1 \in W_1$,
$$[L(-1),\mathcal{Y}(w_1,x)]=\frac{d}{dx}\mathcal{Y}(w_1,x).$$

\item[(iii)](Jacobi identity) For every $w_i \in W_i$, $i=1,2$ and $v \in V$, we have

\bea && x_0^{-1} \delta \left(\frac{x_1-x_2}{x_0}\right)
{Y}(v,x_1) \mathcal{Y}(w_1,x_2)w_2-x_0^{-1} \delta
\left(\frac{-x_2+x_1}{x_0}\right) \mathcal{Y}(w_1,x_2)
{Y}(v,x_1)w_2 \nn &&=x_2^{-1} \delta
\left(\frac{x_1-x_0}{x_2}\right) \mathcal{Y}(Y(v,x_0)w_1,x_2)w_2.
\eea
\end{itemize}
}
\end{definition}
We will denote the vector space of logarithmic intertwining
operators among $W_1$, $W_2$ and $W_3$ by $$I \ {W_3 \choose W_1 \
W_2}.$$ We say that a logarithmic intertwining operator is {\em
genuine} if $(w_1)_\alpha^{(k)}$ is nonzero for some $k \geq 1$
and some $w_1 \in W_1$. A genuine logarithmic intertwining
operator of type ${ W_3 \choose W_1 \ W_2}$ is called {\em hidden}
if two modules $W_i$ and $W_j$ are ordinary $V$-modules, where $\{
i,j \} \subset \{1,2,3\}$.

\begin{definition}
{\em A logarithmic intertwining operator $\mathcal{Y}$ of type $
{W_3 \choose W_1 \ W_2}$ is said to be {\em strong} if \be {\rm
depth}(\mathcal{Y}):=\label{sup} {\rm sup} \{ k : \
{(w_1)}_{\alpha}^{(k)} w_2 \neq 0, \alpha \in \mathbb{C}, w_1 \in
W_1, w_2 \in W_2, k \in \mathbb{N} \} \ee is finite.}
\end{definition}
In other words, for strong logarithmic intertwining operators the
powers of $\lo(x)$ are globally bounded.
\begin{definition} \label{milas-def}
{\em A logarithmic intertwining operator $\mathcal{Y}$  is said to
be {\em locally strong} if
$$\mathcal{Y}(w_1,x) \in {\rm Hom}(W_2,W_3)\{x\}[\lo(x)], \ \ {\rm for} \ \ {\rm every} \ w_1 \in W_1.$$
}
\end{definition}

\noindent The previous definition was used in \cite{M1} as the
definition of logarithmic intertwining operators. In this paper we
shall only study strong intertwining operators. From the previous
definitions we clearly have a chain of embeddings \be \label{hier}
I_{ord} {W_3 \choose W_1 \ W_2} \subseteq I_{st} {W_3 \choose W_1
\ W_2} \subseteq I_{lst} {W_3 \choose W_1 \ W_2} \subseteq I {W_3
\choose W_1 \ W_2 }, \ee where $I_{st}$, $I_{lst}$ and $I_{ord}$
stand for the vector space of strong, locally strong and ordinary
intertwining operators, respectively.

\begin{proposition} \label{strong}
\begin{itemize}
\item[(i)] Every strong logarithmic intertwining operator
$\mathcal{Y}$ of type ${ W_3 \choose W_1 \ W_2}$ defines an
ordinary intertwining operator of the same type. \item[(ii)] Every
locally strong logarithmic intertwining operator $\mathcal{Y}$
among a triple of finitely generated modules is strong.
\end{itemize}
\end{proposition}
\noindent {\em Proof.} Let $k={\rm depth}(\mathcal{Y})$. From
$$\mathcal{Y}(w,x)=\sum_{i=0}^k \mathcal{Y}^{(i)}(w,x) \lo^i(x)$$
it is clear that the truncation condition and Jacobi identity hold
for $\mathcal{Y}^{(k)}$. From $\frac{d}{dx} \lo^k(x)=\frac{k}{x}
\lo^{k-1}(x)$ it follows that
$$[L(-1),\mathcal{Y}^{(k)}(w,x)]=\frac{d}{dx} \mathcal{Y}^{(k)}(w,x).$$

To prove (ii) it suffices to assume that $W_1$ is finitely
generated. Let $\{w_{1,1},...,w_{1,m} \}$ be a generating set of
$W_1$ so that for every $w_1 \in W_1$, there exist $v_i \in V$ and
$n_i \in \mathbb{Z}$ such that $w_1=\sum_{i=1}^k (v_i)_{n_i}
w_{1,i}$. Let
$$\mathcal{Y}(w_{1,i},x)=\sum_{i=1}^{r_i} \mathcal{Y}^{(i)}(w_{1,i},x)
\lo^i(x), \ \ r_i \in \mathbb{N}.$$ The Jacobi identity now gives
\bea \label{15} && \mathcal{Y}(Y(v_{i},x_0) w_{1,i},x) \nn &&
={\rm Res}_{x_1} x_0^{-1} \delta \left(\frac{x_1-x}{x_0}\right)
Y(v_i,x_1)\mathcal{Y}(w_{1,i},x)-x_0^{-1} \delta
\left(\frac{-x+x_1}{x_0}\right) \mathcal{Y}(w_{1,i},x) Y(v_i,x_1).
\eea After we take ${\rm Coeff}_{x_0^{-n_i-1}}$ in (\ref{15}) we
see that the powers of $\lo(x)$ in $\mathcal{Y}((v_i)_{n_i}
w_{1,i},x)$ are bounded by the highest power of $\lo(x)$ in
$\mathcal{Y}(w_{1,i},x)$. Because of the finiteness of the
generating set $W_1$, the powers of $\lo(x)$ in $\mathcal{Y}$ are
globally bounded. \epfv

\begin{remark}
{\em In the previous proposition we considered a canonical map
$$I_{st} {W_3 \choose W_1 \ W_2 } \longrightarrow I_{ord} {W_3 \choose
W_1 \ W_2 },$$ sending a strong logarithmic inertwining operator
to its "top" (or the depth) component. This map is clearly
surjective, but it is far from being injective (see Lemma
\ref{lemma-1}). Therefore, in general, the logarithmic fusion
rules (i.e., ${\rm dim} \ I \ { W_3 \choose W_1 \ W_2 }$) are not
the same as the nonlogarithmic fusion rules (i.e., ${\rm dim} \
I_{ord} \ { W_3 \choose W_1 \ W_2 }$). }
\end{remark}

\begin{remark}
{\em It is tempting to relax the condition in (\ref{map}) and
assume instead that \be \label{mapf} \mathcal{Y}( \cdot, x) \cdot
: \ W_1 \otimes W_2 \longrightarrow W_3 \{x\} [[\lo(x)]]. \ee We
will show (cf. Section 8) that there is a downside for doing that.
One could even put further restriction on (\ref{map}) and force a
freshman calculus "formula" \be \label{freshman} e^{\lo(x)}=x. \ee
This identity, in our formal variable setup, is unnatural and it
should be avoided. The main problem is that $e^{\lo(x)}-x$ is a
unit in the formal ring $\mathbb{C}[[x, \lo(x)]]$.

%(viz. $e^{-\lo(x)} \sum_{n=0}^\infty
%e^{-n \lo(x)} x^n$ is its multiplicative inverse.
}
\end{remark}

\section{An extended Heisenberg Lie algebra}

Let $\goth{h}$ be a finite-dimensional complex abelian Lie algebra
with an inner product $(\cdot,\cdot)$. Consider the central
extension of the affinization of $\goth{h}$, denoted by
$\hat{\goth{h}}$, generated by $h(n):=h \otimes t^n$, $n \in
\mathbb{Z}$, $h \in \goth{h}$, with the bracket relations: \bea
&& [a(m),b(n)]= m(a,b)\delta_{m+n,0}C, \nonumber \\
&& [C,h(m)]=0, \nonumber \eea where $a,b \in \goth{h}$, $m,n \in
\mathbb{Z}$ and $C$ is the central element.

For purposes of this paper we will assume that ${\rm
dim}(\goth{h})=1$ with a fixed unit vector $h$, so that
$$[h(m),h(n)]=m \delta_{m+n,0}C.$$ The Lie algebra $\hat{\goth{h}}$
is an example of an extended Heisenberg Lie algebra (of course,
the Heisenberg Lie algebra associated with $\goth{h}$ does not
involve $h(0)$).
Let
$$\hat{\goth{h}}=\hat{\goth{h}}_{<0} \oplus \hat{\goth{h}}_{>0} \oplus \mathbb{C}h(0) \oplus \mathbb{C}C,$$
where $\hat{\goth{h}}_{<0}$ and $\hat{\goth{h}}_{>0}$ are defined
as usual. Denote by $\mathcal{C}_k$ the category of restricted
$\mathbb{Z}$-graded, $\hat{\goth{h}}$--modules of level $k$, i.e.,
the category of $\mathbb{Z}$-graded modules $W=\coprod_{n \in
\mathbb{Z}} W_n$,
$$h(n) W_m \subseteq W_{m-n}, \ \ {\rm for} \ \ {\rm every} \ \ m,n \in \mathbb{Z},$$
such that there exists $N \in \mathbb{N}$, so that $W_n=0$ for $n
< N$, and the central element $C$ acts as the multiplication with
$k$.

Let us denote by $M(k) \equiv \mathcal{U}(\hat{\goth{h}}_{<0})
\cdot {1}$ essentially unique irreducible lowest weight module of
level $k$, where $h(n) \cdot {\bf 1}=0$, for $n \geq 0$ and the
grading is the obvious one (see \cite{FLM}).

Then we have a version of the Stone-Von Neumann theorem for the
extended Heisenberg algebra $\hat{\goth{h}}$:
\begin{lemma} \label{desheis}
Let $W$ be a restricted $\hat{\goth{h}}$-module of level $k$. Then
\be \label{SVN} W \cong M(k) \otimes \Omega(W), \ee where
$\Omega(W)=\{ w \in W : h(n) w=0, n > 0 \}$ (the vacuum space of
$W$) is $\goth{h}$-stable.
\end{lemma}
\noindent {\em Proof:} It is known (see Theorem 1.7.3 in
\cite{FLM}, for instance) that every restricted module for the
Heisenberg algebra $\tilde{\goth{h}}=\hat{\goth{h}} \ \backslash \
\mathbb{C}h(0)$ admits a decomposition (\ref{SVN}) where
$\Omega(W)$ is the vacuum space of $W$. Now, $h(0)$ commutes with
$h(n)$ for every $n$, hence it preserves the vacuum space
$\Omega(W)$. \epfv

\section{Logarithmic $M(1)_a$-modules}

\noindent The Heisenberg vertex operator algebra is omnipresent in
conformal field theory and representation theory. Despite of its
simplicity, it is the main tool for building (more interesting)
rational vertex operator algebras (e.g., lattice VOAs \cite{D},
\cite{FLM}). Moreover, Heisenberg vertex operator algebras also
contain many interesting subalgebras (e.g., $W$-algebras
\cite{FKRW}).

It is well-known (cf. \cite{FrB}, \cite{FLM}, \cite{K2}) that
$(M(1),Y,{\bf 1}, \omega_a)$ has a vertex operator algebra
structure of central charge $c=1-12a^2$, where
$$\omega_a=\frac{h(-1)^2{\bf 1}}{2}+ah(-2){\bf 1}.$$
This VOA will be denoted by $M(1)_a$. It is known that $M(1)_a$
has infinitely many inequivalent irreducible modules, which can be
easily classified \cite{LL}. For every highest weight irreducible
$M(1)_a$-module $W$ there exists $\lambda \in \mathbb{C}$ such
that $W \cong M(1)_a \otimes \Omega_{\lambda}$, where
$\Omega_\lambda$ is one-dimensional and $h(0)$ acts as the
multiplication with $\lambda$. Such a module will be denoted by
$M(1, \lambda)_a$. The restricted $\hat{\goth{h}}$-modules are
essentially logarithmic $M(1)_a$-modules. Since $M(1)_a$ is at the
same time $\hat{\goth{h}}$, $M(1)_a$ and Virasoro algebra module,
we stress some differences under taking the contragradient module.
If we denote by $M(1,\lambda)^*$ the dual $\hat{\goth{h}}$-module
of $M(1,\lambda)$, under the standard anti-involution $h(n)
\mapsto -h(-n)$, then we have
$$M(1,\lambda)^* \cong M(1,-\lambda),$$
viewed as $\hat{\goth{h}}$. On the other hand, if we denote by
$M(1,\lambda)'_a$ the contragradient $M(1)_a$-module (or the dual
Virasoro algebra module) of $M(1,\lambda)_a$ (cf. \cite{FHL}),
then we have \be \label{cont-iso} M(1,\lambda)'_{a} \cong
M(1,2a-\lambda)_a. \ee
 Let us also mention the automorphism $\tau$ of
$M(1)_a$, uniquely determined by $\tau (h(-n_1) \cdots h(-n_k){\bf
1})=(-1)^k h(-n_1) \cdots h(-n_k){\bf 1}$. This map does not
preserve the Virasoro element $\omega_a$, unless $a=0$.

As in Lemma \ref{desheis} we now describe all logarithmic modules
for the vertex operator algebra $M(1)_a$.

%In what follows we shall
%often use $\Omega$\Omega(W)= \Omega(M(1)_a \otimes { 1} \otimes
%\Omega$.
\begin{proposition}
Suppose that $W$ is a logarithmic $M(1)_a$--module. Then, viewed
as an $\hat{\goth{h}}$-module, $W \cong M(1)_a \otimes \Omega(W)$.
\end{proposition}
\noindent {\em Proof:} Every weak $M(1)_a$--module carries a level
one representation of $\hat{\goth{h}}$ via the expansion
$$Y(h(-1){\bf 1},x)=\sum_{n \in \mathbb{Z}} h(n)x^{-n-1}.$$
Because $W$ is logarithmic, it is also a restricted
$\mathbb{Z}$-graded $\hat{\goth{h}}$-module. An application of
Lemma \ref{desheis} yields $W=M(1)_a \otimes \Omega(W)$, for some
vacuum space $\Omega(W)$.
%Clearly, from the definition of $v_n w=0$ for $v
%\in M(1)_a$, $w \in W$ and $n$ large enough.
%Consider a module $\mathcal{U}(\hat{\goth{h}}) \cdot v$
%Because the $M(1)$-module is restricted it follows that Let $W$ be
%a logarithmic $M$--module (meaning it admits generalized
%eigenspaces decomposition with respect to $L(0)$).
%From
%$$L(0)=\frac{1}{2} \sum_{n \in \mathbb{Z}} \nb h(-n)h(n) \nb-a h(0)$$
%we see that \be \label{ex-lo} L(0)|_{\Omega}=\frac{1}{2}h(0)^2- a
%h(0). \ee If $\Omega$ is a finite-dimensional complex vector space
%then $L(0)|_{\Omega}$ also admits a Jordan form. The proof
%follows.
\epfv

\begin{corollary} \label{hha}
Suppose that $(\frac{h^2(0)}{2}-a h(0))|_{\Omega}$ admits a Jordan
block of size at least $2$. Then $M(1)_a \otimes \Omega$ is a
genuine logarithmic $M(1)_a$-module.
\end{corollary}
%Let $\Lambda=\{ \lambda_1,\lambda_2,\ldots \}$ and $m_i$ be as
%above. Then one can easily obtain the Jordan block decomposition
%with respect to the action of $L(0)$. In particular, if $\lambda_i
%\neq -\lambda_j$ for all pairs $(i,j)$ then the set
%$$\Lambda'=\{\frac{\lambda^2_1}{2},\frac{\lambda^2_2}{2},\ldots \}$$
%and the numbers $m_i$ and partitions $(n_1^{i},n_2^{2},...)$ give
%sizes of all the Jordan blocks (as in the case of $h(0)$).
$ $
In the case of irreducible $M(1)$-modules it possible to classify
all intertwining operators and the corresponding fusion rules. The
following result seems to be known (see \cite{FrB} for instance).
\begin{proposition}
Let $M(1,\lambda)_a$, $M(1,\tau)_a$ and $M(1,\nu)_a$ be three
(ordinary) $M(1)_a$-modules. Then the vector space of ordinary
intertwining operators $I_{ord} \ {M(1,{\nu})_a \choose
M(1,{\lambda})_a \ M(1,{\tau})_a}$ is nontrivial if and only if
$\nu=\lambda+\tau$. If so, $I_{ord} \ {M(1,{\lambda+\tau})_a
\choose M(1,{\lambda})_a \ M(1,{\tau})_a}$ is one-dimensional.
%Moreover, \bea && \mathcal{Y}(h(-i_1) \cdots
%h(-i_k)v_{\lambda},x)= \nn && : \left(\frac{d}{dx}\right)^{i_1}
%h(x) \cdots \left(\frac{d}{dx}\right)^{i_k} h(x)
%\mathcal{Y}(v_{\lambda},x) :. \eea
\end{proposition}
In Section 5 we will substantially generalize this result.

\section{Generalized Logarithmic Intertwining Operators}

Here we show that certain spaces of "logarithmic operators" give
rise to (generalized) logarithmic modules. In particular, this
construction will be useful for classification of logarithmic
intertwining operators among certain triples of logarithmic
$M(1)_a$-modules. We will closely follow Li's work \cite{Li} (see
also \cite{Li1}). The following definition is a logarithmic
version of Definition 6.1.1 in \cite{Li1} (see also \cite{Li}).
\begin{definition} \label{genint}
{\em Let $V$ be a vertex operator algebra and $W_1$ and $W_2$ a
pair of logarithmic $V$--modules. An operator valued formal series
$$\phi(x)=\sum_{\alpha \in \mathbb{C}} \sum_{n \in \mathbb{N}}
\phi_\alpha^{(n)} \lo^n(x) x^{-\alpha-1} \in {\rm Hom}(W_1,W_2)\{x
\}[[\lo(x)]]$$
% $\phi(x) \in {\rm
%Hom}(M_1,M_2)\{x\}[\lo(x)]$
is called a {\em generalized logarithmic intertwining operator} of
generalized weight $h$ if it satisfies the following conditions:
%\begin{itemize}
%\item[(i)] There are finitely many complex numbers $h_1$,..,$h_k$,
%such hat $h_i-h_j \neq \mathbb{Z}$, $i \neq j$ and a positive
%integer $M$ such that for every
%$$\phi(x) \in \sum_{k=1}^{n} \sum_{\alpha=0}^M x^{h_k} \lo^\alpha (x)
%{\rm Hom}(W_1,W_2)[[x,x^{-1}]],$$
%$$\phi(x)u \in  \sum_{k=1}^n \sum_{\alpha=0}^M x^{h_k} \lo^\alpha
%(x)W_2((x)).$$
%$$\phi(x)=\sum_{\alpha \in \mathbb{C}} \phi(\alpha)x^{-\alpha-1} \in \sum_{n \in \mathbb{N}} \lo^n(x)
%{\rm Hom}(W_1,W_2)\{x \},$$ satisfying
\begin{itemize}
\item[(i)] For every $w \in W_1$ we have
$$\phi(x)w \in
W_2 \{x \}[\lo(x)],$$ and for every $i$
$$\phi_{\alpha+m}^{(i)}w=0, \ \ {\rm for} \ \ m \ \ {\rm large \ \ enough}.$$
\item[(ii)] $[L(-1),\phi(x)]=\frac{d}{dx}\phi(x),$ \item[(iii)]
For every $v \in V$, there is a positive integer $n_v$ such that
$$(x_1-x_2)^{n_v}Y(v,x_1)\phi(x_2)=(x_1-x_2)^{n_v}\phi(x_2)Y(v,x_1).$$
\item[(iv)] There exists $k \in \mathbb{N}$ such that
$$\left({\rm ad}_{L(0)}-x\frac{d}{dx}-h \right)^k \phi(x)=0,$$
where ${\rm ad}_{L(0)}\phi(x)=[L(0),\phi(x)]$.

\end{itemize}
}
\end{definition}
We will denote by $G_h(W_1,W_2)$ the vector space of  generalized
logarithmic intertwining operators of generalized weight $h$ and
by $G_{{\rm log}}(W_1,W_2)=\oplus_{h \in \mathbb{C}}G_h(W_1,W_2)$
the vector space of generalized logarithmic intertwining
operators.

\begin{proposition} The vector space $G_{{\rm log}}(W_1,W_2)$ has a generalized logarithmic $V$-module structure.
\end{proposition}
\noindent {\em Proof.} We will closely follow Section 6 in
\cite{Li1} (or Theorem 7.1.6 in \cite{Li}), with appropriate
modifications due to logarithms. As in Definition 6.1.2 in
\cite{Li1} we let
$$Y(u,x_0) \circ \phi(x_2)=
{\rm Res}_{x_1} \left( x_0^{-1} \delta
\left(\frac{x_1-x_2}{x_0}\right) Y(u,x_1) \phi(x_2)-x_0^{-1}
\delta \left(\frac{-x_2+x_1}{x_0}\right) \phi(x_2) Y(u,x_1)
\right).$$ The $L(-1)$-property is then proven as in \cite{Li1}.
Let $a \in V$ be a homogeneous vector of weight ${\rm wt}(a)$ and
$\phi(x)$ a generalized logarithmic intertwining operator of
generalized weight $h$, so that
$$({\rm ad}_{L(0)} -x \frac{d}{dx}- h)^k \phi(x)=0,$$ for some
$k \in \mathbb{N}$. We claim that
$$\left({\rm ad}_{L(0)}-h-{\rm wt}(a)-x_0 \frac{\partial}{\partial x_0}- x_2 \frac{\partial}{\partial
x_2} \right)^k Y(a,x_0) \circ \phi(x_2)=0.$$ By using Lemma 7.1.4
in \cite{Li} we get the identity
$$[L(0),Y(a,x_0) \circ \phi(x_2)]$$
$$=\left({\rm wt}(a)+x_0 \frac{\partial}{\partial x_0}
+{x_2} \frac{\partial}{\partial x_2} \right) Y(a,x_0) \circ
\phi(x_2)+Y(a,x_0) \circ \left({\rm ad}_{L(0)}
\phi(x_2)-x_2\frac{\partial}{\partial x_2} \phi(x_2) \right).$$
Thus
$$\left( {\rm ad}_{L(0)}-{\rm wt}(a)-h- x_0 \frac{\partial}{\partial
x_0}-x_2 \frac{\partial}{\partial x_2} \right) Y(a,x_0) \circ
\phi(x_2)=Y(a,x_0) \circ \left({\rm
ad}_{L(0)}-x_2\frac{\partial}{\partial x_2}-h \right) \phi(x_2)
.$$ From the previous formula and property (iv) in Definition
\ref{genint} the claim now follows. Finally, the Jacobi identity
is the verbatim repetition of the proof of Theorem 6.1.7 in
\cite{Li1}. \epfv

In parallel with the ordinary case, generalized logarithmic
intertwining operators are closely related to logarithmic
intertwining operators.
\begin{lemma}
Let $W$ be a logarithmic $V$-module and $\phi \in {\rm
Hom}_V(W,G_{{\rm log}}(W_1,W_2))$. Define \bea && I_\phi(\cdot,x)
: \ W \longrightarrow {\rm Hom}(W_1,W_2)\{x\}[[\lo(x)]], \nonumber
\\
&& I_\phi(w,x)=\phi(w)(x).\nonumber \eea Then $I_\phi(\cdot,x)$ is
a logarithmic intertwining operator of type $ { W_2 \choose W \
W_1}$.
\end{lemma}
\noindent {\em Proof.} Firstly,
$$I_\phi(w,x) w_1 \in W_2\{x\}[\lo(x)].$$
Now we have to check properties (i)-(iii) in Definition
\ref{def-logint}. Let
$$I_{\phi}(w,x)=\sum_{\alpha \in \mathbb{C}}
\sum_{n \geq 0} w^{(n)}_\alpha \lo(x)^n x^{-\alpha -1}.$$ Then the
truncation property
$$w_{\alpha+m}^{(i)} w_1=0 \ \ {\rm for} \ m >> 0$$
is just a consequence of Definition \ref{genint},(i). Similarly,
the $L(-1)$-property holds. The Jacobi identity is then proven as
in formula (6.2.3) in \cite{Li1}. \epfv

The map between ${\rm Hom}_V(W,G_{\rm log}(W_1,W_2))$ and $I \ {
W_2 \choose W \ W_1}$ defined in the previous lemma is clearly
injective. Conversely, to every logarithmic intertwining operator
$\mathcal{Y}$ of type ${ W_2 \choose W \ W_1}$ we associate a map
$\psi \in {\rm Hom}_V(W,G_{{\rm log}}(W_1,W_2))$, via
$$\psi(w)=\mathcal{Y}(w,x), \  \ w \in W.$$
Combined together we obtain:
\begin{theorem}
Let $W$, $W_1$ and $W_2$ be logarithmic $V$-modules. Then the
vector space $I \ {W_2 \choose W \ W_1}$ is naturally isomorphic
to ${\rm Hom}_V(W,G_{{\rm log}}(W_1,W_2))$.
\end{theorem}
The next result is a logarithmic analogue of Li's Theorem 7.3.1
\cite{Li} (after Tsuchiya and Kanie who proved an important
spacial case $V=L_{sl_2}(k,0)$, $k \in \mathbb{N}$). We do not
need this theorem in full generality so we just focus on a special
case $V=M(1)_a$.
\begin{theorem} \label{thm5}
Let $W_1$ and $W_2$ be two logarithmic $M(1)_a$-modules. Let
$\Omega$ be a finite-dimensional $\goth{h}$-module and
$\mathcal{Y}( \cdot,x)$ a linear map from $\Omega$ to ${\rm
Hom}(W_1,W_2)[[{\rm log}(x)]]\{x\}$ satisfying the truncation
condition, the $L(-1)$-property and \bea \label{jaca} &&
(x_1-x_2)^{{\rm
wt}(a)-1}Y(a,x_1)\mathcal{Y}(w,x_2)-(-x_2+x_1)^{{\rm wt}(a)-1}
\mathcal{Y}(w,x_2) Y(a,x_1)\nn &&=x_1^{-1} \delta \left(
\frac{x_2}{x_1} \right) \mathcal{Y}(a_{{\rm wt}(a)-1} \cdot w,x),
\eea for every homogeneous $a \in V$ and $w \in \Omega$. Then
$\mathcal{Y}$ extends uniquely to an intertwining operator of type
$ \ds{{W_2 \choose M(1)_a \otimes \Omega \ \ W_1}}$.
\end{theorem}
\noindent {\em Proof.} From (\ref{jaca}) it follows that
$\mathcal{Y}(u,x)$ is a generalized intertwining operator for
every $u \in \Omega$.
%As in \cite{Li} we consider
%$$\bar{\Omega}=\{ \mathcal{Y}(a,x) : a \in U \} \subset
%G_{log}(W_1,W_2).$$
Moreover, $\Omega$ is also an $A(M(1)_a)$-module (where the Zhu's
algebra $A(M(1)_a)$ is just the polynomial ring in one variable).
Now we may proceed as in \cite{Li}. The linear map $\mathcal{Y}$
extends to an intertwining map from $M(1)_a \otimes \Omega$ (The
universal Verma $M(1)_a$-module $\bar{M}(\Omega)$ appearing in
Li's theorem is  $M(1)_a \otimes \Omega$.) If there is another
intertwining operator $\mathcal{Y'}$ extending
$\mathcal{Y}|_{\Omega}$, then $\mathcal{Y}-\mathcal{Y}'$ would be
trivial on $\Omega$. But $M(1)_a \otimes \Omega$ is generated by
$\Omega$, so $\mathcal{Y}=\mathcal{Y}'$. \epfv

%Even if $\mathcal{Y}$ is only nontrivial, Theorem \ref{thm5} (with
%small modifications) would yield a logarithmic intertwining
%operator of type ${W_2 \choose M(1)_a \otimes \Omega' \ W_1}$, for
%some $\goth{h}$-submodule $\Omega' \subseteq \Omega$.

\section{Logarithmic intertwining operators among $M(1)_a$-modules}

In this section we will give a sharp upper bound on the dimension of
the vector space of strong logarithmic intertwining operators among
certain logarithmic $M(1)_a$-modules. From now on we shall assume
that every logarithmic $M(1)_a$-module is of the form $M(1)_a
\otimes \Omega$, where $\Omega$ is a finite-dimensional
$\goth{h}$-module such that $(h(0)-\lambda)^n|_{\Omega}=0$ for some
$\lambda$ and $n$ large enough. Then Proposition 1.10 in \cite{M1}
implies that every intertwining operator among a triple of such
modules is strong.  If we remove the finite-dimensionality condition
on $\Omega$ it is not hard to construct logarithmic intertwining
operators that are neither strong nor locally strong.

We prove a few lemmas first.
\begin{lemma}
Suppose that $(h(0)-\lambda)^{m_1}|_{\Omega_1}=0$ and
$(h(0)-\nu)^{m_2}|_{\Omega_2}=0$ for some $m_i \in\mathbb{N}$,
$i=1,2$. Let $ \ds{\mathcal{Y} \in I \ {W \choose M(1)_a \otimes
\Omega_1 \ M(1)_a \otimes \Omega_2}}$.
Then for every $w_1 \in M(1)_a \otimes \Omega_1$ and $w_2 \in
M(1)_a \otimes \Omega_2$ we have
$$(h(0)-\lambda-\nu)^{m_1+m_2-1} \mathcal{Y}(w_1,x)w_2=0.$$
Moreover, if $\Omega_1$ and $\Omega_2$ are one-dimensional
$\goth{h}$-modules, then there are no genuine logarithmic
intertwining operators of type $\ds{{ W \choose M(1)_a \otimes
\Omega_1 \ M(1)_a \otimes \Omega_2}}$.
\end{lemma}
\no {\em Proof.} Let $w_1 \in M(1)_a \otimes \Omega_1$ and $w_2
\in M(1)_a \otimes \Omega_2$ so that $(h(0)-\lambda)^{m_1} \cdot
w_1=(h(0)-\nu)^{m_2} \cdot w_2=0$. From the Jacobi identity it
follows that \bea && (h(0)-\lambda-\nu)^{m_1+m_2-1}
\mathcal{Y}(w_1,x)w_2 \nonumber \\
&& =\sum_{{\tiny \begin{array}{c} i_1 \geq 0, i_2 \geq 0
\\ i_1+i_2=m_1+m_2-1 \end{array}}} {m_1+m_2-1 \choose i_1}
\mathcal{Y}((h(0)-\lambda)^{i_1}w_1,x)(h(0)-\nu)^{i_2}
w_2.\nonumber \eea Now, it is easy to see that every term on the
right hand side is zero.

If $m_1=m_2=1$ then $m_1+m_2-1=1$, so $h(0)$ is diagonalizable on
the image of $\mathcal{Y}$, but so is $L(0)=\frac{1}{2} h(0)^2- a
h(0) + \sum_{n
> 0} h(-n)h(n)$. Now, apply Proposition 1.10 in \cite{M1}. \epfv

Let us recall (cf. \cite{M1}, \cite{HLZ}) that for every $w_i \in
W_i$ of generalized weight $h_i$, $i=1,2$ and every $\mathcal{Y}
\in I \ {W_3 \choose W_1 \ W_2}$, the vector
$$(w_1)^{(k)}_{\alpha} w_2$$
is of generalized weight $h_1+h_2-\alpha-1$ (independently of
$k$).

\begin{lemma} \label{mn}
Suppose that
$(h(0)-\lambda)^{m_1}|_{\Omega_1}=(h(0)-\nu)^{m_2}|_{\Omega_2}=0$
and let \\ $\ds{\mathcal{Y} \in I \ { W \choose M(1)_a \otimes
\Omega_1 \ M(1)_a \otimes \Omega_2}}$. Then
$$\mathcal{Y}(w_1,x)w_2 \in x^{\mu \lambda} W((x)) \oplus x^{\mu \lambda} \lo(x)W((x)) \oplus \cdots \oplus
 x^{\mu \lambda} \lo(x)^{m_1+m_2-2}W((x)),$$
for every $w_i \in M(1) \otimes \Omega_i$, $i=1,2$. Equivalently,
$${\rm depth}(\mathcal{Y}) \leq m_1+m_2-2.$$
\end{lemma}
\noindent {\em Proof.} The proof goes by induction on $m_1+m_2
\geq 2$. For $m_1+m_2=2$ the statement holds by the previous
lemma. From the same lemma and \be \label{ex-lo}
L(0)|_{\Omega_1}=\frac{1}{2}h(0)^2- a h(0)\ee it follows that
$L(0)|_{M(1)_a \otimes \Omega_1}$ (resp. $L(0)|_{M(1)_a \otimes
\Omega_2}$) does not admit a Jordan block of size larger than
$m_1$ (resp. $m_2$). For the induction step we apply the bracket
relation between $L(0)$ and $\mathcal{Y}(w_1,x)$ and an elementary
ODE argument as in Proposition 1.10, \cite{M1}.
\epfv

\begin{lemma} \label{nil-lemma}
Let $\mathcal{Y}$ and $\Omega_i$ be as in Lemma \ref{mn}. Suppose
further that $\lambda=0$. Then for every $w_i \in M(1)_a \otimes
\Omega_i$, $i=1,2$ we have \be \label{nilpotint}
\mathcal{Y}(w_1,x)w_2 \in W ((x)) \oplus \lo(x) W(( x )) \oplus
\cdots \oplus \lo^{m_1-1}(x) W (( x )). \ee If $\nu=0$, then
(\ref{nilpotint}) holds with $m_1$ replaced by $m_2$.
\end{lemma}
\noindent {\em Proof.} By using the isomorphism $ I \ {W \choose
W_1 \ W_2} \cong I \ {W \choose W_2 \ W_1}$ \cite{M1}, it is
sufficient to consider the $\lambda=0$ case, so that
$h(0)|_{\Omega_1}$ is a nilpotent operator. We prove the formula
(\ref{nilpotint}) by induction on $m_1$. Firstly, let $m_1=1$, so
$M(1)_a \otimes \Omega_1 \cong M(1)_a$. In this case we have to
show that there are no genuine logarithmic intertwining operators
(i.e., $\mathcal{Y}(w_1,x)w_2 \in W((x))$).
%Let $\mathcal{Y} \in I \
%\ds{{ W \choose M(1)_a \otimes \Omega_1 \ M(1)_a \otimes
%\Omega_2}}$.
From the $L(-1)$-property and $L(-1){\bf 1}=0$, it
follows that
$$\langle w'_3, \mathcal{Y}(L(-1){\bf 1},x) w_2 \rangle=\frac{d}{dx}
\langle w'_3,\mathcal{Y}({\bf 1},x)w_2 \rangle=0, \ \ w'_3 \in
W'_3.$$ Thus $\mathcal{Y}({\bf 1},x)$ is a constant term
(operator) and it does not involve powers or $\lo(x)$. Similarly,
from the Jacobi identity, it follows that $\mathcal{Y}(w,x)w_1$
does not involve nonzero powers of $\lo(x)$ for every $w \in
M(1)_a$ and $w_1 \in M(1)_a \otimes \Omega_2$.

Now, suppose that (\ref{nilpotint}) holds for every $m < m_1$. For
$w_1 \in \Omega_1$, and $w_2 \in M(1)_a \otimes \Omega_2$ and
$w'_3 \in \Omega(W')=\Omega(W)'$ we clearly have
$$\langle w'_3, \mathcal{Y}(L(-1)w_1,x)w_2 \rangle=\frac{d}{dx} \langle w'_3,\mathcal{Y}(w_1,x)w_2 \rangle.$$
On the other hand, the Jacobi identity gives
$$\mathcal{Y}(L(-1)w_1,x)=\mathcal{Y}(h(-1)h(0)w_1,x)=\nb h(x)
\mathcal{Y}(h(0) w_1,x) \nb,$$ so that
$$ \frac{d}{dx} \langle w'_3,\mathcal{Y}(w_1,x)w_2 \rangle=\langle w'_3,\mathcal{Y}(h(0)w_1,x) h(0) x^{-1} w_2 \rangle $$
or \be \label{ind-log} x\frac{d}{dx} \langle
w'_3,\mathcal{Y}(w_1,x)w_2 \rangle=\langle
w'_3,\mathcal{Y}(h(0)w_1,x) h(0) w_2 \rangle. \ee If we denote
$$\bar{\Omega}(W_1)=\{h(0) w_1 : \ w_1 \in \Omega_1 \},$$
then $h(0)^{m_1-1}|_{\bar{\Omega}(W_1)}=0$ so by induction
hypothesis the right hand side in (\ref{ind-log}) is of the form
$P({\rm log}(x))$, where ${\rm deg}(P) \leq m_1-2$. Therefore
$\langle w'_3,\mathcal{Y}(w_1,x)w_2 \rangle$ is a polynomial in
$\lo(x)$ of degree at most $m_1-1$. Now, the Jacobi identity
implies that in $\langle w'_3, \mathcal{Y}(w_1,x)w_2 \rangle$ the
powers of ${\rm log}(x)$ are bounded by $m_1-1$ for every $w_1 \in
M(1)_a \otimes \Omega_1$. The proof follows. \epfv

\begin{lemma} \label{lemma-1}
Let $\mathcal{Y} \in I \ {W_3 \choose W_1 \ W_2 }$ be a strong
logarithmic intertwining operator of depth $k$ among an arbitrary
triple of logarithmic $V$-modules. Then
$$\mathcal{Y}_{-1}( \cdot,x):=\sum_{i=0}^{k-1} (i+1)
\mathcal{Y}^{(i+1)}( \cdot,x) \lo^i(x),$$ defines a strong
intertwining logarithmic operator of depth $k-1$.
\end{lemma}
\noindent {\em Proof.} The truncation condition and the Jacobi
identity clearly hold for $\mathcal{Y}_{-1}( \cdot,x)$. It remains
to prove the $L(-1)$-property. If we distribute the
$L(-1)$-property for $\mathcal{Y}(\cdot,x)$ among
$\mathcal{Y}^{(i)}( \cdot,x)$ we get
$$[L(-1),\mathcal{Y}^{(k)}(w,x)]=\frac{d}{dx}
\mathcal{Y}^{(k)}(w,x),$$ and
$$[L(-1),\mathcal{Y}^{(i)}(w,x)]=\frac{d}{dx}
\mathcal{Y}^{(i)}(w,x)+\frac{i+1}{x} \mathcal{Y}^{(i+1)}(w,x),$$
for $i \leq k-1$. Thus \bea &&
[L(-1),\mathcal{Y}_{-1}(w,x)]=\sum_{i=0}^{k-1}
(i+1)[L(-1),\mathcal{Y}^{(i+1)}(w,x)] \lo^i(x) \nn && =k
\frac{d}{dx} \left( \mathcal{Y}^{(k)}(w,x) \right) \lo^{k-1}(x)+
\sum_{i=0}^{k-2} (i+1) \left(\frac{d}{dx}
\mathcal{Y}^{(i+1)}(w,x)+\frac{i+2}{x} \mathcal{Y}^{(i+2)}(w,x)
\right) \lo^i(x) \nn && =\sum_{i=0}^{k-1} (i+1) \frac{d}{dx}
\left(\mathcal{Y}^{(i+1)}(w,x) \lo^{i}(x) \right). \nonumber \eea
\epfv

The following result gives a sharp upper bound on the depth of
(strong) logarithmic intertwining operators among a triple of
logarithmic $M(1)_a$-modules.
\begin{theorem}  \label{uppertheorem} Let $h(0)|_{\Omega_1}$ and $h(0)|_{\Omega_2}$ such that
$$(h(0)-\lambda)^{m_1}|_{\Omega_1}=(h(0)-\nu)^{m_2}|_{\Omega_2}=0,$$
and
$$(h(0)-\lambda)^{m_1-1}|_{\Omega_1} \neq 0,
\ (h(0)-\nu)^{m_2-1}|_{\Omega_2} \neq 0,$$ for some $\lambda, \nu
\in \mathbb{C}$. Then we have
\begin{itemize}
\item[(i)] For every $\mathcal{Y} \in I \ \ds{{W \choose M(1)_a
\otimes \Omega_1 \ \ M(1)_a \otimes \Omega_2}}$, \bea
\label{interval} 0 \leq {\rm depth}(\mathcal{Y}) \leq k=\left\{
\begin{array}{ccc} m_1+m_2-2 & {\rm for} & \lambda \nu \neq 0, \\
m_1-1 & {\rm for} & \lambda=0 \ \ {\rm and} \ \ \nu \neq 0 \\
m_2-1 &  {\rm for} & \lambda \neq 0 \ \ {\rm and} \ \ \nu=0
\\ {\rm min}(m_1-1,m_2-1) & {\rm for} & \lambda=\nu=0 \end{array} \right.
\eea \item[(ii)]There exists a canonical embedding \bea
\label{canonical}
 I \ {W \choose M(1)_a
\otimes \Omega_1 \ \ M(1)_a \otimes \Omega_2} \hookrightarrow {\rm
Hom}_{\goth{h}}(\Omega_1 \otimes \Omega_2,\Omega(W))^{\oplus^{k
}}, \eea where $k$ is as in (\ref{interval}).
% \eea
%
% =\left\{
%begin{array}{ccc} m_1+m_2-2 & {\rm for} & \lambda \nu \neq 0, \\
%m_1-1 & {\rm for} & \lambda=0 \ \ {\rm and} \ \ \nu \neq 0 \\
%m_2-1 &  {\rm for} & \lambda \neq 0 \ \ {\rm and} \ \ \nu=0
%\\ {\rm min}(m_1,m_2) & {\rm for} & \lambda=\nu=0 \end{array} \right.
%\eea
\item[(iii)] The range of ${\rm depth}(\mathcal{Y})$ in
(\ref{interval}) is the best possible. More precisely, for every
nonnegative integer $m \leq k$ there exists a logarithmic
intertwining operator of depth exactly $m$.
\end{itemize}
\end{theorem}
\noindent {\em Proof.} Here we prove (i) and (ii) only. We will
complete the proof of (iii) in Section 7.

The assertion (i) follows from Lemma \ref{mn} and \ref{nil-lemma}.
Let $\mathcal{Y} \in I \ \ds{{W \choose M(1)_a \otimes \Omega_1 \ \
M(1)_a \otimes \Omega_2}}$. Then $\mathcal{Y}$ admits a canonical
expansion
$$\mathcal{Y}(w_1,x)=\sum_{i=0}^{m_1+m_2-2} \mathcal{Y}^{(i)}(w_1,x) \lo^i(x),$$
for every $w_1 \in M(1)_a \otimes \Omega_1$. By Lemma \ref{mn},
for every $i$,
$$\mathcal{Y}^{(i)}(w_1,x)w_2 \in x^{\lambda \nu} W((x)).$$
Also, for $w_1$ and $w_2$ satisfying
$(L(0)-h_1)^{m_1}w_1=(L(0)-h_2)^{m_2}w_2=0,$ the vector $(w
_1)^{(i)}_\alpha w_2$ is homogeneous of generalized weight
$h_1+h_2-\alpha-1$ for every $i$. Suppose that $w_1 \in \Omega_1$
and $w_2 \in \Omega_2$ are of generalized weight
$\frac{\lambda^2}{2}-a \lambda$ and $\frac{\nu^2}{2}- a \nu$,
respectively. Then for every $i$,
$$w_1 \otimes w_2 \mapsto (w_1)^{(i)}_{-\lambda \nu-1} (w_2).$$
defines an $\goth{h}$-module map
$$ {F}^{(i)}_{\mathcal{Y}} : \Omega_1 \otimes \Omega_2 \longrightarrow \Omega(W).$$
We claim that \bea && F \ : \ I \ { W \choose M(1)_a \otimes
\Omega_1 \ \ M(1)_a \otimes \Omega_2} \longrightarrow {\rm
Hom}_{\goth{h}}(\Omega_1 \otimes \Omega_2,\Omega(W))^{\oplus^k}
\nn &&
F(\mathcal{Y})=(F^{(0)}_{\mathcal{Y}},...,F^{(m_1+m_2-2)}_{\mathcal{Y}}),
\nonumber \eea defines an embedding. Suppose that
$F^{(i)}_{\mathcal{Y}}=0$ for every $i$. We have to show that
$\mathcal{Y} \equiv 0$.

Firstly, we prove that $\mathcal{Y}(w_1,x)w_2=0$ for every $w_1 \in
\Omega_1$ and $w_2 \in \Omega_2$. In order to prove that we observe
first that the contragradient (or dual) module of $W$, denoted by
$W'$, is again a logarithmic module generated by
$\Omega(W')=\Omega(W)'$. By the assumption we have
$$\langle w'_3,\mathcal{Y}(w_1,x)w_2 \rangle=0, \ \ w_3 \in
\Omega(W').$$

Furthermore, for every $w'_3 \in \Omega(W')$ and $n \geq 1$ we
have
$$\langle h(-n) \cdot w'_3,\mathcal{Y}(w_1,x)w_2
\rangle=- \langle w'_3,x^{-n}\mathcal{Y}(h(0) w_1,x)w_2
\rangle=0.$$ Since $W'$ is generated by $h(n)$, $n \leq -1$ from
$\Omega(W')$, the previous formula gives
$$\langle w'_3 , \mathcal{Y}(w_1,x)w_2 \rangle=0,$$
for every $w'_3 \in W'$. Hence, $\mathcal{Y}(w_1,x)w_2=0$ for
every $w_1 \in \Omega_1$ and $w_2 \in \Omega_2$. Finally, $M(1)_a
\otimes \Omega_i$ is generated by $\Omega_i$, so
$\mathcal{Y}(w_1,x)w_2=0$ for every $w_i \in M(1)_a \otimes
\Omega_i$, $i=1,2$.
%
%
%Thus $\mathcal{Y}(w_1,x)w_2=0$. But then $\mathca
%
%
%We will show that then $(w_1)^{(i)}_{\alpha}=0$ for every $w_1 \in
%\Omega_1$, $\alpha$ and $i$, which implies $\mathcal{Y} \equiv 0$.
%Suppose that there exist $w_1 \in \Omega_1$, $w_2 \in M(1)_a
%\otimes \Omega_2$, $i$ and $\alpha$ such that
%$$(w_1)^{(i)}_\alpha w_2 \neq 0,$$
%and $(w_1)^{(i)}_{\alpha+j}w_2=0,$ for $j>0$. Choose $w_2$ to be
%of the highest possible generalized weight. Then
%$(w_1)^{(i)}_\alpha w_2$ is a highest weight vector for the
%Heisenberg algebra. Thus $(w_1)^{(i)}_\alpha w_2$ is of
%generalized weight $\lambda \nu=\frac{\lambda^2}{2}-a \lambda+{\rm
%wt}(w_2)-\alpha-1$. But
%$$(w_1)^{(i)}_\alpha h(n)w_20$ for all
%$n>1$, implies that $(w_1)^{(i)}_\alpha $
% so $\alpha=-\lambda \nu -1$ (the weight of
%$w_1$ is fixed). But this is contradiction with
%$F_{\mathcal{Y}}^{(i)}=0$.
% It is not hard to see that
%then that $\alpha=-\lambda \nu -1$. But we assumed that
%$F_{\mathcal{Y}}^{(i)}=0$ for every $i$, so for every $w_2 \in
%\Omega_2$ and $w_1 \in \Omega_1$ we have
%$\mathcal{Y}(w_1,x)w_2=0$. But then because of
%$[h(n),\mathcal{Y}(w_1,x)]=x^n \mathcal{Y}(h(0) \cdot w_1,x)$ it
%follows that $\mathcal{Y}(w_1,x)w_2=0$ for every $w_1 \in M(1)_a
%\otimes \Omega_1$ and $w_2 \in M(1)_a \otimes \Omega_2$.
This
proves the injectivity.
% and
%$$(w_1)^{(i)}_{\alpha+j} w_2=0$$ for all $j>0$. Now for $w_2 \in
%\Omega_2$ from the Jacobi identity it follows
%$$h(n) (w_1)^{(i)}_\alpha w_2=0, \ \ {\rm for} \ n \geq 1,$$
%o that $w=(w_1)^{(i)}_\alpha w_2=0$ for every
%
%
%is a highest weight vector for $\hat{\goth{h}}$. Thus $w$ is a
%vector in $\Omega(W)$. Hence the kernel of $F$ is trivial. This
%proves (ii).

In Section 7 we will complete the proof of (iii). Now, let us
assume that there exists $\mathcal{Y}$ and $W$ such that ${\rm
depth}(\mathcal{Y})=k$ where $k$ is as in (\ref{interval}) (i.e.,
${\rm depth}(\mathcal{Y})$ reaches its upper bound). Then to
construct $\mathcal{Y}_i$ with ${\rm depth}(\mathcal{Y}_i)=k-i$,
for $i=1,...,k$ we simply apply Lemma \ref{lemma-1}. \epfv

%Suppose now that we were given a $(m_1+m_2-1)$-tuple of
%$\goth{h}$-equivariant maps $(F^{(0)},...,F^{(m_1+m_2-2)})$.
%Suppose that $F^{(m_1+m_2-2)}$ is nonzero first. Then we can
%construct a logarithmic intertwining operator
%$\mathcal{Y}_1(\cdot,x)$ of depth $m_1+m_2-2$ such that
%$$F_{\mathcal{Y}_1}^{(m_1+m_2-2)}=F^{(m_1+m_2-2)}.$$
%%$$w_1 \otimes w_2 \longrightarrow (w_1)^{(m_1+m_2-2)}_{-\lambda \nu -1}(w_2)$$
%%equals the map $F^{(m_1+m_2)}$.
%Consider the difference \be \label{difference}
%{F}^{(m_1+m_2-3)}-{F}^{(m_1+m_2-3)}_{\mathcal{Y}}. \ee Now we
%choose $\bar{\mathcal{Y}}_1$ of depth $m_1+m_2-2$ such that
%$$F_{\bar{\mathcal{Y}}_1}^{(m_1+m_2-2)}=\frac{1}{m_1+m_2-2}
%\left({F}^{(m_1+m_2-3)}-{F}^{(m_1+m_2-3)}_{\mathcal{Y}}\right).$$
%Apply now Lemma \ref{lemma-1} for $\bar{\mathcal{Y}}_1$ and we
%obtain a logarithmic intertwining operator
%$$\mathcal{Y}_{1,2}(\cdot,x):=\mathcal{Y}_1(\cdot,x)+(\bar{\mathcal{Y}}_1)_{-1}(\cdot,x),$$
%which satisfy
%$$F_{\mathcal{Y}_{1,2}}^{(m_1+m_2-2)}=F^{(m_1+m_2-2)}$$
%and
%$$F_{\mathcal{Y}_{1,2}}^{(m_1+m_2-3)}=F^{(m_1+m_2-3)}.$$
%Now, inductively we can construct an intertwining operator of
%$\mathcal{Y}$ such that
%$$F_{\mathcal{Y}}^{(i)}=F^{(i)}, \ \ i=0,...,m_1+m_2-2.$$

%In case $F^{(m_1+m_2-2-j)}=0$ for some $j=0,1,..,k$, we apply
%Lemma \ref{lemma-1} $k$ times and proceed as before. \epfv

\section{The operator "$x^{h(0)}$"}

Let us recall that every endomorphism of a finite-dimensional
complex vector space admits a unique decomposition
$$h(0)=h_s(0)+h_n(0),$$
where $h_s(0)$ and $h_n(0)$ are the semisimple and nilpotent part
of $h(0)$, respectively such that $h_s(0)$ and $h_n(0)$ commute.

The following elementary fact will be of use in the next section:
Let $\Omega$ be a finite dimensional vector space and $h(0)$ an
endomorphism of $\Omega$. Then \be \label{xho} e^{ \lo(x) h_n(0)}
x^{h_s(0)}, \ee is a solution of the ODE
$$x\frac{d}{dx} A(x)=h(0) A(x),$$
where
$$A(x) : \ V \longrightarrow V \{x \}[\lo(x)],$$
is a linear map. The operator valued expression (\ref{xho}) is a
replacement for $x^{h(0)}$ that (in the case when $h(0)$ is
semisimple) appears on many places in the literature. Clearly, the
"operator" $x^{h(0)}$ for $h(0)$ nonsemisimple is not well-defined.

\section{Logarithmic intertwining operators for $M(1)_a$: The proof of
the existence}

\noindent In this section to every
$T_{\Omega_1,\Omega_2}^{\Omega_3} \in {\rm
Hom}_{\goth{h}}(\Omega_1 \otimes \Omega_2, \Omega_3)$ we associate
a linear map
$$\mathcal{Y} ( \cdot, x) \cdot : \Omega_1 \otimes W_2
\longrightarrow W_3\{x\} [\lo(x)],$$
$$W_i=M(1)_a \otimes \Omega_i, \ \ i=2,3,$$
with the following properties: \bea \label{asss} &&
[h(n),\mathcal{Y}(w,x)]=x^{n} \mathcal{Y}( h(0) \cdot w,x),  \ n \in
\mathbb{Z}, \\ \label{ass} && [L(-1),\mathcal{Y}(w,x)]=\frac{d}{dx}
\mathcal{Y}(w,x),  \eea for every $w \in \Omega_1$. We will
construct this map explicitly. Our construction mimics the
well-known formulas when $\Omega_i$ are all one-dimensional, which
will be a special case of our construction.

Firstly, we fix a basis $\mathcal{B}=\{w_1,...,w_n \}$ for
$\Omega_1$ in which $h(0)|_{\Omega_1}$ admits a Jordan form
consisting of a single Jordan block so that \bea \label{jordan}
h(0) \cdot w_1&=& \lambda w_1, \nn h(0) \cdot w_i&=&\lambda
w_i+w_{i-1}, \ {\rm for} \ 2 \leq i \leq n. \eea
Let \bea
\int^{+} h(x)&=&h(0) \lo(x) + \sum_{m
>0} \frac{h(m) x^{-m}}{-m}, \nn \int^{-} h(x)&=&\sum_{m < 0}
\frac{h(m) x^{-m}}{-m}. \nonumber \eea Formal differentiation
yields
$$\frac{d}{dx} \left( \int^+ h(x) + \int^- h(x) \right)=\sum_{n \in \mathbb{Z}} h(n) x^{-n-1}.$$
Let \bea  E^+(\lambda,x)&=& {\rm exp} \left(\sum_{m
>0} \frac{\lambda h(m) x^{-m}}{-m} \right), \nn
E^-(\lambda,x)&=& {\rm exp} \left(\sum_{m < 0} \frac{ \lambda h(m)
x^{-m}}{-m} \right). \nonumber \eea Here we use a slightly
different notation compared with \cite{LL}; $E^\pm(\lambda,x)$ are
usually denoted by $E^\pm(-\lambda,x)$. The following lemma is
easy to prove (see for instance \cite{FLM}, \cite{K}, \cite{LL},
etc.).
\begin{lemma}
\bea \label{l-1} && [L(-1),E^-(\lambda,x)]=(\sum_{n \leq -2}
\lambda h(n)x^{-n-1}) E^-(\lambda,x), \nn &&
[L(-1),E^+(\lambda,x)]= E^+(\lambda,x) (\sum_{n \geq 0} \lambda
h(n) x^{-n-1}). \eea
\end{lemma}

Similarly, \bea \label{EE} && [h(n),E^-(\lambda,x)E^+(\lambda,x)]=
x^{n} \lambda E^-(\lambda,x)E^+(\lambda,x), \ \ {\rm for} \ n \neq
0, \nn && [h(0),E^-(\lambda,x)E^+(\lambda,x)]= 0. \nonumber \eea
Also,
$$[h(0),T_{\Omega_2}^{\Omega_3}(w)]=T_{\Omega_2}^{\Omega_3}( h(0)
\cdot w), \ \ w \in \Omega_1,$$ so that \bea \label{honon} &&
[h(0),E^-(\lambda,x)E^+(\lambda,x)T_{\Omega_2}^{\Omega_3}(w)e^{\lo(x)
\lambda h_n(0)}x^{\lambda h_s(0)}]\nn &&
=E^-(\lambda,x)E^+(\lambda,x)T_{\Omega_2}^{\Omega_3}(h(0) \cdot
w)e^{\lo(x) \lambda h_n(0)}x^{\lambda h_s(0)}. \eea
\begin{lemma} \label{l1lemma}
For every $ w \in \Omega_1$ we have \bea &&
[L(-1),E^-(\lambda,x)E^+(\lambda,x)T_{\Omega_2}^{\Omega_3}(w)e^{\lo(x)
\lambda h_n(0)}x^{\lambda h_s(0)}] \nn &&=\left(\sum_{n \leq -2}
\lambda h(n)x^{-n-1} \right)E^-(\lambda,x) E^+(\lambda,x)
T_{\Omega_2}^{\Omega_3}(w)e^{\lo(x) \lambda h_n(0)}x^{\lambda
h_s(0)} \nn && + E^-(\lambda,x) E^+(\lambda,x)
T_{\Omega_2}^{\Omega_3}(w) e^{\lo(x) \lambda h_n(0)} x^{\lambda
h_s(0)} \left(\sum_{n \geq 0} \lambda h(n) x^{-n-1} \right) \nn &&
+ h(-1)E^-(\lambda,x)E^+(\lambda,x) T_{\Omega_2}^{\Omega_3}(h(0)
\cdot w) e^{\lo(x)\lambda h_n(0)} x^{\lambda h_s(0)}. \eea In
particular if $h(0) \cdot w=\lambda w$, then \bea \label{semisimp}
&&
[L(-1),E^-(\lambda,x)E^+(\lambda,x)T_{\Omega_2}^{\Omega_3}(w)e^{\lo(x)\lambda
h_n(0)}x^{ \lambda h_s(0)}] = \nn && \frac{d}{dx}
\left(E^-(\lambda,x)E^+(\lambda,x)T_{\Omega_2}^{\Omega_3}(w)e^{\lo(x)\lambda
h_n(0)}x^{\lambda h_s(0)}\right). \eea
\end{lemma}
\noindent {\em Proof.} Let us recall that $L(-1)=\frac{1}{2}
\sum_{n \in \mathbb{Z}} \nb h(-n-1) h(n) \nb$, which does not
depend on $a$. For simplicity let
$$A(\lambda,w,x)=E^-(\lambda,x)E^+(\lambda,x)T_{\Omega_2}^{\Omega_3}(w)e^{\lo(x)\lambda
h_n(0)}x^{ \lambda h_s(0)}.$$ By using
$$[L(-1),T_{\Omega_2}^{\Omega_3}(w)]=h(-1)
T_{\Omega_2}^{\Omega_3}( h(0) \cdot w)$$ we have  \bea &&
[L(-1),A(\lambda,w,x)]=\left(\sum_{n \leq -2} \lambda h(n)x^{-n-1}
\right)E^-(\lambda,x) E^+(\lambda,x)
T_{\Omega_2}^{\Omega_3}(w)e^{\lo(x) \lambda h_n(0)}x^{\lambda
h_s(0)} \nn &&+ E^-(\lambda,x) E^+(\lambda,x) \left(\sum_{n \geq
0} \lambda h(n) x^{-n-1} \right) T_{\Omega_2}^{\Omega_3}(w)
e^{\lo(x)\lambda h_n(0)} x^{\lambda h_s(0)} \nn &&  +
E^-(\lambda,x) E^+(\lambda,x) h(-1) T_{\Omega_2}^{\Omega_3}(h(0)
\cdot w) e^{\lo(x)\lambda h_n(0)} x^{\lambda h_s(0)} \nn &&=
\left(\sum_{n \leq -2} \lambda h(n)x^{-n-1} \right)E^-(\lambda,x)
E^+(\lambda,x) T_{\Omega_2}^{\Omega_3}(w)e^{\lo(x) \lambda
h_n(0)}x^{\lambda h_s(0)} \nn && + E^-(\lambda,x) E^+(\lambda,x)
T_{\Omega_2}^{\Omega_3}(w) e^{\lo(x)\lambda h_n(0)} x^{\lambda
h_s(0)} \left(\sum_{n \geq 0} \lambda h(n) x^{-n-1} \right) \nn
&&+ \lambda x^{-1} E^-(\lambda,x) E^+(\lambda,x)
T_{\Omega_2}^{\Omega_3}(h(0) \cdot w) e^{\lo(x)\lambda h_n(0)}
x^{\lambda h_s(0)} \nn
%&& + E^-(\lambda,x) E^+(\lambda,x)
%T_{\Omega_2}^{\Omega_3}(h(0) \cdot w) e^{\lo(x)\lambda h_n(0)}
%x^{\lambda h_s(0)} \left(\sum_{n \geq 0} \lambda h(n) x^{-n-1}
%\right) \nn
&& -x^{-1} \lambda E^-(\lambda,x) E^+(\lambda,x)
T_{\Omega_2}^{\Omega_3}(h(0) \cdot w) e^{\lo(x)\lambda h_n(0)}
x^{\lambda h_s(0)}+ \nn && h(-1)E^-(\lambda,x) E^+(\lambda,x)
T_{\Omega_2}^{\Omega_3}(h(0) \cdot w) e^{\lo(x)\lambda h_n(0)}
x^{\lambda h_s(0)} \nn && =\left(\sum_{n \leq -2} \lambda
h(n)x^{-n-1} \right)E^-(\lambda,x) E^+(\lambda,x)
T_{\Omega_2}^{\Omega_3}(w)e^{\lo(x) \lambda h_n(0)}x^{\lambda
h_s(0)} \nn && + E^-(\lambda,x) E^+(\lambda,x)
T_{\Omega_2}^{\Omega_3}(w) e^{\lo(x) \lambda h_n(0)} x^{ \lambda
h_s(0)} \left(\sum_{n \geq 0} \lambda h(n) x^{-n-1} \right) \nn
&&+ h(-1)E^-(\lambda,x)E^+(\lambda,x) T_{\Omega_2}^{\Omega_3}(h(0)
\cdot w) e^{\lo(x) \lambda h_n(0)} x^{\lambda h_s(0)}. \nonumber
\eea \epfv

\begin{lemma} \label{h0lemma}
$$[h(n),E^-(\lambda,x)E^+(\lambda,x)T_{\Omega_2}^{\Omega_3}(w)e^{\lo(x)\lambda h_n(0)}x^{\lambda h_s(0)}]
$$
$$= \lambda x^{n}
E^-(\lambda,x)E^+(\lambda,x)T_{\Omega_2}^{\Omega_3}(w)e^{\lo(x)
\lambda h_n(0)}x^{\lambda h_s(0)}, \ \ {\rm for} \ n \neq 0,$$
$$[h(0),E^-(\lambda,x)E^+(\lambda,x)T_{\Omega_2}^{\Omega_3}(w)e^{\lo(x) \lambda h_n(0)}x^{\lambda h_s(0)}]
$$
$$=E^-(\lambda,x)E^+(\lambda,x)T_{\Omega_2}^{\Omega_3}(h(0) \cdot w)e^{\lo(x) \lambda h_n(0)}x^{\lambda h_s(0)}.$$
\end{lemma}

The following theorem gives a solution to equations (\ref{asss})
and (\ref{ass}) for a Jordan block of length $n$.

\begin{theorem} \label{jordanthm}
Assume that  $w_i$ satisfy (\ref{jordan}). Let \bea \label{intone}
\mathcal{Y}(w_1,x)=E^-(\lambda,x)E^+(\lambda,x)T_{\Omega_2}^{\Omega_3}(w_1)e^{\lo(x)
\lambda h_n(0)}x^{\lambda h_s(0)}, \eea  and for $2 \leq i \leq
n$, let \bea \label{inti} && \\ && \mathcal{Y}(w_i,x)=\sum_{l=1}^i
\left( \sum_{j=0}^{i-l} \frac{(\int^- h(x))^j}{j!} E^-(\lambda ,x)
E^+(\lambda,x) T_{\Omega_2}^{\Omega_3}(w_l) e^{\lo(x) \lambda
h_n(0)}x^{\lambda h_s(0)} \frac{(\int^+
h(x))^{i-l-j}}{(i-l-j)!}\right). \nonumber \eea Then \bea
\label{inth0} && [L(-1),\mathcal{Y}(w_i,x)]=\frac{d}{dx}
\mathcal{Y}(w_i,x), \nn && [h(n),\mathcal{Y}(w_i,x)]= x^n
\mathcal{Y}(h(0) \cdot w_i,x). \eea
\end{theorem}
\noindent {\em Proof.} We already proved the formula in the $i=1$
case. So we may assume $i \geq 2$. For every $ j \geq 1$ we have
\bea && [L(-1), \frac{(\int^- h(x))^j}{j!} E^-(\lambda ,x)
E^+(\lambda,x) T_{\Omega_2}^{\Omega_3}(w_1) e^{\lo(x) \lambda
h_n(0)}x^{\lambda h_s(0)} \frac{(\int^+
h(x))^{i-l-j}}{(i-l-j)!}]=\nn && \left(\sum_{n \leq -2}
h(n)x^{-n-1} \right) \frac{(\int^- h(x))^{j-1}}{(j-1)!}
E^-(\lambda ,x) E^+(\lambda,x) T_{\Omega_2}^{\Omega_3}(w_1)
e^{\lo(x) \lambda h_n(0)}x^{\lambda h_s(0)} \frac{(\int^+
h(x))^{i-l-j}}{(i-l-j)!} \nn &&+\frac{(\int^- h(x))^{j}}{(j)!}
\frac{d}{dx} \left(E^-(\lambda ,x) E^+(\lambda,x)
T_{\Omega_2}^{\Omega_3}(w_1) e^{\lo(x) \lambda h_n(0)}x^{\lambda
h_s(0)} \right) \frac{(\int^+ h(x))^{i-l-j}}{(i-l-j)!} \nn && +
\frac{(\int^- h(x))^{j}}{(j)!} \frac{d}{dx} \left(E^-(\lambda ,x)
E^+(\lambda,x) T_{\Omega_2}^{\Omega_3}(w_1) e^{\lo(x) \lambda
h_n(0)}x^{\lambda h_s(0)} \right) \frac{(\int^+
h(x))^{i-l-j-1}}{(i-l-j-1)!} \left( \sum_{n \geq 0} h(n)x^{-n-1}
\right) \nonumber \eea \bea && =\frac{d}{dx} \left( \frac{(\int^-
h(x))^j}{j!} E^-(\lambda ,x) E^+(\lambda,x)
T_{\Omega_2}^{\Omega_3}(w_1) e^{\lo(x) \lambda h_n(0)}x^{\lambda
h_s(0)} \frac{(\int^+ h(x))^{i-l-j}}{(i-l-j)!} \right) \nn &&-
h(-1) \left( \frac{(\int^- h(x))^{j-1}}{(j-1)!} E^-(\lambda ,x)
E^+(\lambda,x) T_{\Omega_2}^{\Omega_3}(w_1) e^{\lo(x) \lambda
h_n(0)}x^{\lambda h_s(0)} \frac{(\int^+ h(x))^{i-l-j}}{(i-l-j)!}
\right).  \nonumber \eea

For $2 \leq l \leq i-1$, by using Lemma \ref{l1lemma} we have \bea
&& {[L(-1), \frac{(\int^- h(x))^j}{j!} E^-(\lambda ,x)
E^+(\lambda,x) T_{\Omega_2}^{\Omega_3}(w_l) e^{\lo(x) \lambda
h_n(0)}x^{\lambda h_s(0)} \frac{(\int^+ h(x))^{i-l-j}}{(i-l-j)!}]}
\nn &&=\left(\sum_{n \leq -2} h(n)x^{-n-1} \right) \frac{(\int^-
h(x))^{j-1}}{(j-1)!} E^-(\lambda ,x) E^+(\lambda,x)
T_{\Omega_2}^{\Omega_3}(w_l) e^{\lo(x) \lambda h_n(0)}x^{\lambda
h_s(0)} \frac{(\int^+ h(x))^{i-l-j}}{(i-l-j)!} \nn &&
+\frac{(\int^- h(x))^{j}}{j!} \biggl\{ \left(\sum_{n \leq -2}
\lambda h(n)x^{-n-1} \right)E^-(\lambda,x) E^+(\lambda,x)
T_{\Omega_2}^{\Omega_3}(w_l)e^{\lo(x) \lambda h_n(0)}x^{\lambda
h_s(0)} \nn && + E^-(\lambda,x) E^+(\lambda,x)
T_{\Omega_2}^{\Omega_3}(w_l) e^{\lo(x) \lambda h_n(0)} x^{\lambda
h_s(0)} \left(\sum_{n \geq 0} \lambda h(n) x^{-n-1} \right) \nn &&
+ h(-1)E^-(\lambda,x)E^+(\lambda,x) T_{\Omega_2}^{\Omega_3}(h(0)
\cdot w_l) e^{\lo(x)\lambda h_n(0)} x^{\lambda h_s(0)} \biggr\}
\frac{(\int^+ h(x))^{i-l-j}}{(i-l-j)!} \nn && + \frac{(\int^-
h(x))^{j}}{(j)!} \left(E^-(\lambda ,x) E^+(\lambda,x)
T_{\Omega_2}^{\Omega_3}(w_l) e^{\lo(x) \lambda h_n(0)}x^{\lambda
h_s(0)} \right) \frac{(\int^+ h(x))^{i-l-j-1}}{(i-l-j-1)!} \left(
\sum_{n \geq 0} h(n)x^{-n-1} \right) \nonumber \eea \bea &&=
\left(\sum_{n \leq -2} h(n)x^{-n-1} \right) \frac{(\int^-
h(x))^{j-1}}{(j-1)!} E^-(\lambda ,x) E^+(\lambda,x)
T_{\Omega_2}^{\Omega_3}(w_l) e^{\lo(x) \lambda h_n(0)}x^{\lambda
h_s(0)} \frac{(\int^+ h(x))^{i-l-j}}{(i-l-j)!} \nn &&
+\frac{(\int^- h(x))^{j}}{j!} \biggl\{ \left(\sum_{n \leq -1}
\lambda h(n)x^{-n-1} \right)E^-(\lambda,x) E^+(\lambda,x)
T_{\Omega_2}^{\Omega_3}(w_l)e^{\lo(x) \lambda h_n(0)}x^{\lambda
h_s(0)} \nn && + E^-(\lambda,x) E^+(\lambda,x)
T_{\Omega_2}^{\Omega_3}(w) e^{\lo(x) \lambda h_n(0)} x^{\lambda
h_s(0)} \left(\sum_{n \geq 0} \lambda h(n) x^{-n-1} \right) \nn &&
+ h(-1)E^-(\lambda,x)E^+(\lambda,x)
T_{\Omega_2}^{\Omega_3}(w_{l-1}) e^{\lo(x)\lambda h_n(0)}
x^{\lambda h_s(0)} \biggr\} \frac{(\int^+ h(x))^{i-l-j}}{(i-l-j)!}
\nn && + \frac{(\int^- h(x))^{j}}{(j)!} \left(E^-(\lambda ,x)
E^+(\lambda,x) T_{\Omega_2}^{\Omega_3}(w_l) e^{\lo(x) \lambda
h_n(0)}x^{\lambda h_s(0)} \right) \frac{(\int^+
h(x))^{i-l-j-1}}{(i-l-j-1)!} \left( \sum_{n \geq 0} h(n)x^{-n-1}
\right) \nn && =\frac{d}{dx} \left( \frac{(\int^- h(x))^j}{j!}
E^-(\lambda ,x) E^+(\lambda,x) T_{\Omega_2}^{\Omega_3}(w_l)
e^{\lo(x) \lambda h_n(0)}x^{\lambda h_s(0)} \frac{(\int^+
h(x))^{i-l-j}}{(i-l-j)!} \right) \nn &&-h(-1)\frac{(\int^-
h(x))^{j-1}}{(j-1)!} E^-(\lambda ,x) E^+(\lambda,x)
T_{\Omega_2}^{\Omega_3}(w_l) e^{\lo(x) \lambda h_n(0)}x^{\lambda
h_s(0)} \frac{(\int^+ h(x))^{i-l-j}}{(i-l-j)!} \nn && +
h(-1)\frac{(\int^- h(x))^{j}}{j!} E^-(\lambda,x)E^+(\lambda,x)
T_{\Omega_2}^{\Omega_3}(w_{l-1}) e^{\lo(x)\lambda h_n(0)}
x^{\lambda h_s(0)} \frac{(\int^+
h(x))^{i-l-j}}{(i-l-j)!}.\nonumber \eea

Also, by Lemma \ref{l1lemma} we have \bea && {[L(-1), E^-(\lambda
,x) E^+(\lambda,x) T_{\Omega_2}^{\Omega_3}(w_i) e^{\lo(x) \lambda
h_n(0)}x^{\lambda h_s(0)}]} \nn && =\frac{d}{dx} \left(
E^-(\lambda ,x) E^+(\lambda,x) T_{\Omega_2}^{\Omega_3}(w_i)
e^{\lo(x) \lambda h_n(0)}x^{\lambda h_s(0)} \right) \nn &&+h(-1)
E^-(\lambda,x)E^+(\lambda,x) T_{\Omega_2}^{\Omega_3}(w_{i-1})
e^{\lo(x)\lambda h_n(0)}.\nonumber \eea

By combining the previous three formulas we obtain \bea &&
\sum_{l=1}^i \sum_{j=0}^{i-l} {[L(-1), \frac{(\int^- h(x))^j}{j!}
E^-(\lambda ,x) E^+(\lambda,x) T_{\Omega_2}^{\Omega_3}(w_l)
e^{\lo(x) \lambda h_n(0)}x^{\lambda h_s(0)} \frac{(\int^+
h(x))^{i-l-j}}{(i-l-j)!}]} \nn && =\sum_{l=1}^i \sum_{j=0}^{i-l}
\frac{d}{dx} \left(\frac{(\int^- h(x))^j}{j!} E^-(\lambda ,x)
E^+(\lambda,x) T_{\Omega_2}^{\Omega_3}(w_l) e^{\lo(x) \lambda
h_n(0)}x^{\lambda h_s(0)} \frac{(\int^+
h(x))^{i-l-j}}{(i-l-j)!}\right)+ \nn && -h(-1) \sum_{l=1}^{i-1}
\sum_{j=1}^{i-l} \left( \frac{(\int^- h(x))^{j-1}}{(j-1)!}
E^-(\lambda ,x) E^+(\lambda,x) T_{\Omega_2}^{\Omega_3}(w_l)
e^{\lo(x) \lambda h_n(0)}x^{\lambda h_s(0)} \frac{(\int^+
h(x))^{i-l-j}}{(i-l-j)!} \right) \nn && + h(-1) \sum_{l=2}^i
\sum_{j=0}^{i-l} \left( \frac{(\int^- h(x))^j}{j!} E^-(\lambda ,x)
E^+(\lambda,x) T_{\Omega_2}^{\Omega_3}(w_{l-1}) e^{\lo(x) \lambda
h_n(0)}x^{\lambda h_s(0)} \frac{(\int^+ h(x))^{i-l-j}}{(i-l-j)!}
\right) \nonumber \eea \bea && =\frac{d}{dx} \mathcal{Y}(w_i,x)
\nn && -h(-1) \sum_{l=1}^{i-1} \sum_{j=1}^{i-l} \left(
\frac{(\int^- h(x))^{j-1}}{(j-1)!} E^-(\lambda ,x) E^+(\lambda,x)
T_{\Omega_2}^{\Omega_3}(w_l) e^{\lo(x) \lambda h_n(0)}x^{\lambda
h_s(0)} \frac{(\int^+ h(x))^{i-l-j}}{(i-l-j)!} \right) \nn && +
h(-1) \sum_{l=2}^i \sum_{j=0}^{i-l} \left( \frac{(\int^-
h(x))^j}{j!} E^-(\lambda ,x) E^+(\lambda,x)
T_{\Omega_2}^{\Omega_3}(w_{l-1}) e^{\lo(x) \lambda
h_n(0)}x^{\lambda h_s(0)} \frac{(\int^+ h(x))^{i-l-j}}{(i-l-j)!}
\right) \nn && =\frac{d}{dx} \mathcal{Y}(w_i,x). \nonumber \eea

Assume now that $n >0$. The formula (\ref{inth0}) certainly holds
for $i=1$. Then for $i \geq 2$, by Lemma \ref{h0lemma}, we get
\bea && [h(n), \sum_{l=1}^i \left( \sum_{j=0}^{i-l} \frac{(\int^-
h(x))^j}{j!} E^-(\lambda ,x) E^+(\lambda,x)
T_{\Omega_2}^{\Omega_3}(w_l) e^{\lo(x) \lambda h_n(0)}x^{\lambda
h_s(0)} \frac{(\int^+ h(x))^{i-l-j}}{(i-l-j)!}\right)] \nn
&&=\sum_{l=1}^{i-1} \left( \sum_{j=1}^{i-l} x^n \frac{(\int^-
h(x))^{j-1}}{(j-1)!} E^-(\lambda,x) E^+(\lambda,x)
T_{\Omega_2}^{\Omega_3}(w_l) e^{\lo(x) \lambda h_n(0)}x^{\lambda
h_s(0)} \frac{(\int^+ h(x))^{i-l-j}}{(i-l-j)!}\right) \nn && +
\sum_{l=1}^i \left( \sum_{j=0}^{i-l} \lambda x^n \frac{(\int^-
h(x))^j}{j!} E^-(\lambda ,x) E^+(\lambda,x)
T_{\Omega_2}^{\Omega_3}(w_l) e^{\lo(x) \lambda h_n(0)}x^{\lambda
h_s(0)} \frac{(\int^+ h(x))^{i-l-j}}{(i-l-j)!}\right) \nn &&
=\sum_{l=1}^{i-1} \left( \sum_{j=0}^{i-l-1} x^n \frac{(\int^-
h(x))^{j}}{j!} E^-(\lambda.x) E^+(\lambda,x)
T_{\Omega_2}^{\Omega_3}(w_l) e^{\lo(x) \lambda h_n(0)}x^{\lambda
h_s(0)} \frac{(\int^+ h(x))^{i-l-j-1}}{(i-l-j-1)!}\right) \nn && +
\sum_{l=1}^i \left( \sum_{j=0}^{i-l} \lambda x^n \frac{(\int^-
h(x))^j}{j!} E^-(\lambda ,x) E^+(\lambda,x)
T_{\Omega_2}^{\Omega_3}(w_l) e^{\lo(x) \lambda h_n(0)}x^{\lambda
h_s(0)} \frac{(\int^+ h(x))^{i-l-j}}{(i-l-j)!}\right) \nn && = x^n
\mathcal{Y}(w_{i-1},x)+ x^n \lambda \mathcal{Y}(w_i,x). \nonumber
\eea
%\sum_{l=1}^i \left( \sum_{j=0}^{i-l} \frac{(\int^- h(x))^j}{j!}
%E^-(\lambda ,x) E^+(\lambda,x) T_{\Omega_2}^{\Omega_3}(h(0) \cdot
%w_l) e^{\lo(x) \lambda h_n(0)}x^{\lambda h_s(0)} \frac{(\int^+xs
%h(x))^{i-l-j}}{(i-l-j)!}\right). \eea
Consequently,
$$[h(n),\mathcal{Y}(w_i,x)]=x^n \mathcal{Y}(h(0) \cdot
w_i,x).$$ It is easy to see that the previous formula holds for
$n<0$ as well. Finally, by Lemma \ref{h0lemma}
$$[h(0),\mathcal{Y}(w_i,x)]=\mathcal{Y}(h(0) \cdot w_i,x).$$
\epfv

\begin{remark}
{\em Even though the intertwining operator  $\mathcal{Y}(w,x)$
associated with $T_{\Omega_1,\Omega_2}^{\Omega_3}$ was defined for
a single Jordan block only, it is now straightforward to define
$\mathcal{Y}(w,x)$ for an arbitrary finite-dimensional
$\goth{h}$-module $\Omega_1$.}
\end{remark}

Finally, we have a description of logarithmic intertwining
operators among a triple of logarithmic modules with
finite-dimensional vacuum spaces.
\begin{theorem} \label{pre-upper}
Let $W_1=M(1)_a \otimes \Omega_1$, $W_2=M(1)_a \otimes \Omega_2$
and $W_3 =M(1)_a \otimes \Omega_3$ as above and
$$T_{\Omega_1,\Omega_2}^{\Omega_3} \in {\rm
Hom}_{\goth{h}}(\Omega_1,{\rm Hom}(\Omega_2,\Omega_3)),$$ then
$\mathcal{Y}(w,x)$ associated with
$T_{\Omega_1,\Omega_2}^{\Omega_3}$ as in Theorem \ref{jordanthm},
and defined for $w \in \Omega_1$ only, extends uniquely to a
logarithmic intertwining operator of type ${ W_3 \choose W_1 \ W_2
}$.
%Moreover, every logarithmic intertwining
%operator of ${ W_3 \choose W_1 \ W_2 }$ arise in this way.
\end{theorem}
\noindent {\em Proof.} Because of Theorem \ref{jordanthm} and
Theorem \ref{thm5}, it only remains to show (\ref{jaca}). The
formula (\ref{jaca}) clearly holds in the case when $v=h(-1){\bf
1}$ by Theorem \ref{jordanthm}. But the vertex operator algebra is
generated by the vector $h(-1){\bf 1}$ so (\ref{jaca}) holds for
every homogeneous $v \in M(1)_a$ (cf. Proposition 7.1.5 and 7.3.3
in \cite{Li} ). \epfv

\noindent {\em Proof of Theorem \ref{uppertheorem} (iii).} We
choose $\Omega(W)=\Omega_1 \otimes \Omega_2$, such that
$T_{\Omega_1, \Omega_2}^{\Omega_1 \otimes \Omega_2}=Id \in {\rm
Hom}_{\goth{h}} (\Omega_1 \otimes \Omega_2,\Omega_1 \otimes
\Omega_2)$. Let $\mathcal{Y} \in I \ds{ { M(1)_a \otimes \Omega_1
\otimes \Omega_2 \choose M(1)_a \otimes \Omega_1 \ \ M(1)_a
\otimes \Omega_2}}$ as constructed in Theorem \ref{pre-upper} and
Theorem \ref{jordanthm}. It is not hard to see that for such
$\mathcal{Y}$ and $\lambda \nu \neq 0$, we have ${\rm
depth}(\mathcal{Y})=m_1+m_2-2$, because of (\ref{intone}) and
(\ref{inti}). Similarly for $\lambda=0$ or $\nu=0$. \epfv

\section{Mock logarithmic intertwiners among ordinary modules}

In this section (which is completely independent from the rest of
the paper) we construct certain operators related to logarithmic
intertwining operators studied in earlier sections. These operators
involve logarithms but operate among ordinary $M(1)_a$-modules. The
next result has already been proven.
\begin{proposition} The only logarithmic intertwining
operators among $M(1,\lambda)_a$, $M(1,\nu)_a$ and $M(1,
\lambda+\nu)_a$ are the ordinary intertwining operators.
\end{proposition}
%\noindent {\em Proof:} This is simply a consequence of the
%Proposition 1.10 in \cite{M1}.
%\epfv

\begin{proposition}
Suppose that in Definition \ref{def-logint}  we allow
$\mathcal{Y}( \cdot ,x) \cdot$ to satisfy (\ref{mapf}). Under this
new definition, there exists a nontrivial (mock) logarithmic
intertwining operator among every triple of ordinary
$M(1)_a$-modules $M(1,\lambda)_a$, $M(1,\nu)_a$ and
$M(1,\lambda+\nu)_a$, provided that $\lambda \nu \neq 0$.
\end{proposition}
\no {\em Proof.} Let $\mathcal{Y} \in I \ {M(1,\lambda+\nu)_a
\choose M(1,\lambda)_a \ M(1,\nu)_a}$. It is easy to see that the
operator
$$\bar{\mathcal{Y}}( \cdot ,x)=\mathcal{Y}( \cdot ,x)x^{-
\lambda h(0)}$$ satisfies the Jacobi identity, but it doesn't
satisfy the $L(-1)$-property. For the same reason
$$\mathcal{Y}(\cdot ,x)x^{- \lambda h(0)}h(0)^k \lo(x)^k $$
satisfies the Jacobi identity (notice that $h(0)^k$ acting on
$M(1,\nu)_a$ is merely $\nu^k$), but again it does not satisfy the
$L(-1)$-property. Finally, we consider
$$\mathcal{Y}_{log}( \cdot ,x)= \mathcal{Y}( \cdot ,x)x^{- \lambda h(0)} e^{ \lambda h(0)
\lo(x)},$$ where $$e^{\lambda h(0) \lo(x)}=\sum_{n=0}^\infty
\frac{\lambda^n h(0)^n \lo(x)^n}{n!}.$$ Now, as before
$\mathcal{Y}_{log}( \cdot, x)$ satisfies the Jacobi identity but
also the $L(-1)$-property
$$[L(-1),\mathcal{Y}_{log}(w,x)]=\frac{d}{dx}
\mathcal{Y}_{log}(w,x),$$ which follows from $L(-1)$-property for
$\mathcal{Y}( \cdot, x)$ and the formula
$$\frac{d}{dx}(x^{-\lambda h(0)} e^{\lambda h(0) \lo(x)})=0.$$
\epfv

\section{Indecomposable and logarithmic representations of the Virasoro algebra of central charge $c=1$}

In Section 7 we described logarithmic intertwining operators
associated with logarithmic $M(1)_a$-modules. Here we restrict our
construction to an important vertex operator subalgebra
$L(1,0)=\mathcal{U}(Vir) \cdot {\bf 1} \subset M(1)_0$, where
$L(c,h)$, $(c,h) \in \mathbb{C}^2$ denote the irreducible lowest
weight irreducible module for the Virasoro algebra of central
charge $c$ and lowest conformal weight $h$ \cite{KR}. For
simplicity, we shall also use $M(1)$ instead of $M(1)_0$. Results
from this section form logarithmic extension of several results
from \cite{M2}.

Let us introduce some notation. Let $W$ be a $Vir$-module. By
${\bullet}$ we will denote a lowest weight vector inside $W$ such
that $U(Vir_{\leq 0}) \cdot {\bullet}$ is an irreducible lowest
weight module (every $\bullet$ is of course a singular vector in
$W$). By a $\diamond$ we will denote a vector that becomes a lowest
weight vector in the quotient of $W$ by moding out the submodule
generated by all lowest weight vectors. These vectors will be called
subsingular vectors. Similarly $\rhd$ will denote a vector that
becomes a lowest vector after quotienting with the submodule of $W$
generated by all lowest weight vectors and all subsingular vectors.
Such a vector is called a sub-subsingular vector. One can continue
in this manner and introduce vectors that become lowest weight
vectors after quotienting with the submodule generated by the lowest
weight, subsingular and sub-subsingular vectors, but we shall not
need those in the paper. An arrow $\diamond \rightarrow \bullet$
indicates that $\bullet$ is contained in the submodule generated by
$\diamond$, etc. We also recall here that cosingular vectors are
those vectors that are being mapped to singular vectors in the
contragradient module. There is also a pairing between singular
vectors and equivalence classes of cosingular vectors in a module.

Let us illustrate these definition with an example: Let $W$ be a
$Vir$-module such that \be \label{cc0} 0 \rightarrow L(1,m^2)
\rightarrow W \rightarrow L(1,(m+1)^2) \rightarrow 0, \ee is a
nonsplit extension, so that $W$ is generated by a subsingular
vector of weight $(m+1)^2$. This extension may be visualized as
follows
$$
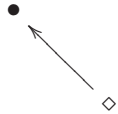
\xymatrix{ {\bullet} & \\ & \ar[ul] {\diamond}}$$
where the arrow pointing up indicates that the conformal weight of
$\bullet$ is smaller than of $\diamond$. More complicated diagrams
will appear later.

In the previous example we implicitly assumed the following result
(see \cite{M2}).
\begin{proposition} \label{milas1}
For every $k,m \in \mathbb{Z}$, we have
$${\rm Ext}^1_{Vir,L(0)}(L(1,k^2),L(1,m^2))=\mathbb{C},$$
if and only if $|k-m|=1$. In all other cases ${\rm
Ext}^1_{Vir,L(0)}$ is trivial.
\end{proposition}
%The previous theorem does not apply for logarithmic modules.
Our aim is to determine Virasoro submodule structure of $M(1)
\otimes \Omega$ (viewed as a ${ Vir}$-module) for some special
$h(0)|_{\Omega}$.

The first part in the following result is well-known (see for
instance \cite{KR}). For the second part see \cite{DG} or
\cite{M2}.
\begin{theorem} \label{griess}
Viewed as a Virasoro module $M(1)$ decomposes as a direct sum of
irreducible Virasoro modules \be \label{decomp} M(1)=\bigoplus_{m
= 0}^\infty L(1,m^2). \ee If $u^m$ denotes the lowest weight
vector (unique up to a nonzero scalar) of $L(1,m^2)$, then the
Virasoro module generated by $Y(u^n,x)u^m$ decomposes as
$$L(1,(m+n)^2) \oplus L(1,(m+n-2)^2) \oplus \cdots \oplus
L(1,(m-n)^2).$$
\end{theorem}
The previous theorem can be used for construction of some
intertwining operators among irreducible $L(1,0)$-modules. This
construction relies on non-vanishing of certain $3j$-symbols
\cite{DG}, \cite{M2}. For instance, $ {\rm dim} \ I \ { L(1,1)
\choose L(1,1) \ L(1,1) }=1$, but this "fusion rule" is not
covered by Theorem \ref{griess}.

Here is a useful consequence of Theorem \ref{griess}
\begin{corollary} \label{lemma-mm1}
The Virasoro module generated by
$$\{ h(-n)u^{m}: n \in \mathbb{Z} \}$$
is isomorphic to \be \label{mm1} L(1,(m-1)^2) \oplus L(1,(m+1)^2),
\ee for $m \geq 1$. If $m=0$, the first summand in (\ref{mm1}) is
trivial.
\end{corollary}
\noindent {\em Proof.} We just have to observe that $h(-1){\bf
1}=u^1$, so that $Y(u^1,x)=\sum_{n \in \mathbb{Z}} h(n) x^{-n-1}.$
Now, apply the previous theorem. \epfv

%\begin{lemma} \label{lemma-mm2}
%Let $W$ be a non-split, graded Virasoro extension of $L(1,m^2)$ by
%$L(1,(m+1)^2)$. Then $W$ generated by a cosingular vector of
%weight $m^2$.
%\end{lemma}
%\noindent {\em Proof.} Let $w_{1,m^2} \in W$ be any nonzero vector
%of weight $m^2$ (cosingular vector in $W$). We also denote by
%$w_{1,(m+1)^2} \in W$ the lowest weight vector of $L(1,(m+1)^2)
%\subset W$. It suffices to show that there exists $a \in U(Vir)$
%such that $a \cdot w \in L(1,(m+1)^2)$. Suppose that this is not
%true. Choose $a_{sing} \in U(Vir_{<0})$ which corresponds to the
%singular vector in the Verma module $M(1,m^2)$ of weight $(m+1)^2$
%(such a vector certainly exists \cite{KR}). If $ a_{sing} \cdot
%w_{1,m^2}=0$, then $U(Vir) \cdot w_{1,m^2} \cong
%M(1,m^2)/I_{max}$, where $I_{max}$ is the maximal submodule of
%$M(1,m^2)$. Thus, $U(Vir) \cdot w \cong L(1,m^2)$ and $(U(Vir)
%\cdot w) \cap L(1,(m+1)^2)=0$. But then $W = L(1,m^2) \oplus
%L(1,(m+1)^2)$. Contradiction. \epfv
%
%It is not hard to see that every nontrivial extension of $W$ (as
%in the previous lemma) by $L(1,(m-1)^2)$ is also a cyclic module,
%and its only proper submodules are $L(1,(m-1)^2)$, $L(1,(m+1)^2)$,
%and their direct sum. This extension can be visualized as a wedge
%(cf. Section 10).

\begin{lemma} \label{duality}
Let $h(0)|_{\Omega}$ be a nilpotent operator. Then $M(1) \otimes
\Omega$ is a self-dual $M(1)$-module, i.e., $(M(1) \otimes
\Omega)' \cong M(1) \otimes \Omega$. Clearly, the same is true if
$M(1) \otimes \Omega$ is viewed as a Virasoro algebra module.
\end{lemma}
\noindent {\em Proof.} Every logarithmic $M(1)$-module $W$ is
uniquely determined by the $h(0)$-action on $\Omega(W)$. Thus, two
logarithmic $M(1)$-modules $M(1) \otimes \Omega_1$ and $M(1)
\otimes \Omega_2$ are equivalent if and only if there exists $\Psi
: \Omega_1 \longrightarrow \Omega_2$ such that
$h(0)|_{\Omega_1}=\Psi^{-1} h(0)|_{\Omega_2} \Psi$. The module
$(M(1) \otimes \Omega)'$ is isomorphic to $M(1) \otimes \Omega'$,
where the action of $h(0)$ on the dual space $\Omega'$ is given by
$$\langle h(0) \cdot w', w \rangle=- \langle w',h(0) \cdot w \rangle, \ w \in \Omega, w' \in \Omega',$$
so that $h(0)|_{\Omega'}=-h^*(0)$, where $h^*(0)$ is the dual map.
The operator $h(0)$ is nilpotent, thus $-h^*(0)$ and consequently
$h(0)|_{\Omega'}$ are nilpotent as well. But $-h(0)^*$ and $h(0)$
admit the same Jordan form, so there exists $\Psi$ with wanted
properties.
%
%The key point is to use the nilpotency of $h(0)$: The operator
%$h(0)^T|_{\Omega}$ is similar to $h(0)|_{\Omega}$, but because of
%nilpotency $-h(0)^T|_{\Omega}$ is also similar to $h(0)_{\Omega}$.
\epfv

\begin{theorem} \label{2step}
Let $\Omega$ be a two-dimensional space and $h(0)|_{\Omega} \neq
0$, $h(0)^2 |_{\Omega}= 0$. Then viewed as a Virasoro module,
$M(1) \otimes \Omega$ is generated by a sequence of subsingular
\footnote{These subsingular vectors are also cosingular.}vectors
as on the following diagram \be \label{embed} \xymatrix{ &
{\bullet} \ar@{.>}[r] &
{} \ar[ld] \diamond \\
{} & {\bullet} \ar@{.>}[r] & {\diamond} \ar[lu] \ar[ld]
\\ {} & {\bullet} \ar@{.>}[r] & {\diamond} \ar[lu] \ar[ld]  \\ {} &
{\bullet} \ar@{.>}[r]
& {} \ar[lu] \ar[ld] \diamond  \\
{} & {\bullet} \ar@{.>}[r] & {\diamond} \ar[lu] \ar[ld]
\\ {} & .. . &  ...\ar[lu] } \ee
where the $s$-th $\bullet$ and $\diamond$, counting from the top,
have conformal weight $(s-1)^2$, $s \geq 1$. Here dotted arrows
indicate the action of the transpose of $h(0)$, which uniquely
determines every $\diamond$.
\end{theorem}
\noindent {\em Proof.} Firstly, we may and will choose a basis $\{
w_1,w_2 \}$ of $\Omega$ such that $h(0) \cdot w_1=0$, $h(0) \cdot
w_2=w_1$, so that $h(0)^T \cdot w_1=w_2$ and $h(0)^T \cdot w_2=0$.
Since $h(0)^2|_{\Omega}=0$, the module $M(1) \otimes \Omega$ is
$L(0)$-diagonalizable. From the Virasoro algebra embedding $M(1)
\hookrightarrow M(1) \otimes \Omega$, where we use the
identification $M(1)=\mathcal{U}(\hat{\goth{h}}) \cdot w_1$, and
Theorem \ref{griess}, it is clear that $M(1) \otimes \Omega$
contains a sequence of singular vectors of weight $m^2$ for every
$m \geq 0$. These singular vectors, displayed in the left column
of (\ref{embed}) by $\bullet$ are determined up to a constant. Let
$u^m=P_m(h) w_1$ denotes such a vector of weight $m^2$, where
$P_m(h)$ is a polynomial in $h(-i)$, $i \geq 1$, of degree $m^2$.
As we already mentioned there are also vectors in $M(1) \otimes
\Omega$ that become singular after quotienting with $M(1)$. These
vectors are uniquely determined if we assume that every $\diamond$
is obtained from $\bullet$ by applying $h(0)^T$ to $u^m$. These
vectors will be denoted by $u^{2,m}$, $m \geq 0$, so that
$u^{2,m}=h(0)^T u^m=P_m(h) \cdot w_2$. It is clear that $M(1)
\otimes \Omega$ is generated by $S=\{ u^m: m \geq 0 \} \cup
\{u^{2,m}: m \geq 0 \}$. It remains to prove that we can reduce
the generating set $S$ down to $\{u^{2,m}: m \geq 0 \}$.
%Clearly, every $u^m$ These vectors are also
%of weight $m^2$.
The short exact sequence of ${\rm Vir}$-modules
$$0 \longrightarrow M(1) \longrightarrow M(1) \otimes \Omega
\stackrel{\pi}{\longrightarrow} M(1) \longrightarrow 0,$$ together
with Theorem \ref{griess} gives
$$0 \longrightarrow \oplus_{m \geq 0} L(1,m^2) \longrightarrow M(1) \otimes \Omega
\stackrel{\pi}{\longrightarrow} \oplus_{m \geq 0} L(1,m^2)
\longrightarrow 0.$$ Thus $M(1) \otimes \Omega$, which is
$L(0)$-diagonalizable, gives a nonzero element in \be
\label{ext-prod} {\rm Ext}^1_{Vir,L(0)}(\oplus_{m \geq 0}
L(1,m^2),\oplus_{n \geq 0} L(1,n^2)) \cong \prod_{m} \prod_{n}
{\rm Ext}^1_{Vir, L(0)}(L(1,m^2),L(1,n^2)). \ee Now, Proposition
\ref{milas1} implies that \be \label{ext-prod1} {\rm
Ext}^1_{Vir,L(0)}(\oplus_{m \geq 0} L(1,m^2),\oplus_{n \geq 0}
L(1,n^2)) \cong \prod_{|m-n|=1} {\rm
Ext}^1_{Vir,L(0)}(L(1,m^2),L(1,n^2)). \ee Already from the
previous formula it is clear that there could be at most two
arrows exiting from $u^{2,m}$.
%Suppose that the image of $\diamond$ in $M(1) \otimes \Omega$ is
%represented by $\bullet$,
%Now, $\pi(M)$ is a submodule of $M(1)$, thus a sum of some
%$L(1,m^2)$. $Ker(\pi(M))$ is also a
%By (\ref{ext-prod1}) there are only two possibilities for arrows
%starting from $u^{2,m}$ of weight $m^2$: towards $u^{m-1}$ or
%towards $u^{m+1}$.
Now, we determine these arrows for every $m$. For $m=0$, there is
precisely one outgoing arrow from $u^{2,0}$ pointing to $u^{1}$.
This follows from $L(-1)w_2=h(-1)h(0)w_2=h(-1)w_1=u^{1}$, where
$u^{1}$ is the
lowest weight vector of $L(1,1) \subset M(1)$.  \\
\no {\em Claim:} For $m \geq 1$, from each $u^{2,m}$ there are
precisely two outgoing arrows; one pointing to $u^{m-1}$ and the
other pointing to $u^{m+1}$.

%First we show that there are arrows pointing up for every $m \geq
%1$. Consider the vector $L(n) u^{2,m}$, for $n \geq 0$. We claim
%that there exists $n>0$ such that $L(n)u^{2,m} \neq 0$.
To see that we write the generator $L(n)$, $n \in \mathbb{Z}$ as
$$h(0)h(n)+\frac{1}{2} \sum_{k+l=n, kl \neq 0} \nb h(k)h(l)
\nb.$$ Now, $$L(n)u^{2,m}=L(n) P_m(h) w_2$$
$$=\left(h(0)h(n)+\frac{1}{2} \sum_{k+l=n, kl \neq 0} \nb h(k)h(l)
\nb \right)P_m(h)w_2=h(n) P_m(h)w_1+\bar{L}(n) P_m(h)w_2,$$ where
$$\bar{L}(n)=\frac{1}{2} \left( \sum_{k+l=n, kl \neq 0} \nb
h(k)h(l) \nb \right).$$ By using Corollary \ref{lemma-mm1}, we
have
$$h(n) P_m(h) w_1 \in L(1,(m-1)^2) \oplus L(1,(m+1)^2), \ \ n \in \mathbb{Z}.$$
Combined with $\bar{L}(n) u^{2,m}=0$, $n \geq 1$, it follows that
there exists $a^+ \in \mathcal{U}(Vir_{> 0})$, such that
$$a^+ \cdot u^{2,m}=u^{m-1}.$$
This proves that there is an arrow pointing to $u^{m-1}$. Let us
recall that under taking the dual the singular vectors in $M(1)
\otimes \Omega$ are mapped to cosingular vectors and vice-versa.
In addition, orientations of arrows are reversed. Now, Lemma
\ref{duality} yields an isomorphism between $M(1) \otimes \Omega$
and its dual, which maps singular vectors $u^{m}$ to cosingular
vectors $w^{2,m}$ and cosingular vectors $u^{2,m}$ to  singular
vectors $w^{m}$ such that $w^{m}$ and $w^{2,m}$ form a Jordan
block with respect to $h(0)$ (i.e., $h(0) \cdot w^{m}=0$, $h(0)
\cdot w^{2,m}=w^{m}$). Clearly, $w^{m}$ generates the lowest
weight module $L(1,m^2)$. Thus, there will be an arrow pointing
down from $w^{2,m-1}$ to $w^{m}$. Now, by using the same argument
as before we argue that there is an arrow pointing up from every
$w^{2,m}$ to $w^{m-1}$.
The proof follows.
%(A different proof of the claim can be given if we employ Lemma
%\ref{duality}.)
\epfv

%Since there is no canonical choice for $u^{2,m}$, it is not true
%that every submodule of $M(1) \otimes \Omega$ is generated by a
%subset of $S$ (e.g., consider the Virasoro submodule generated
%$u^0+u^{2,0}$).

Here is a useful consequence of the previous theorem. For
simplicity we only consider $\Omega$ with ${\rm dim}(\Omega)=3$.
\begin{corollary} \label{om-3}
Let $\Omega$ be a three dimensional $\goth{h}$-module such that
$h(0)^3|_{\Omega}=0$ and $h(0)^2|_{\Omega} \neq 0.$ Then, viewed
as a $Vir$-module, $M(1) \otimes \Omega$ is generated by a
sequence of sub-subsingular vectors (denoted by $\rhd$) as on the
following diagram \be 
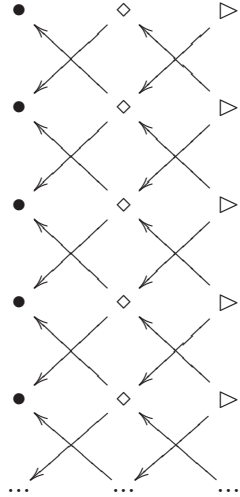
\label{embed2} \xymatrix{{} & {\bullet} &
{\diamond} \ar[ld] & {\rhd}  \ar[ld] \\
{} & {\bullet}  & {\diamond} \ar[lu]\ar[ld] & {\rhd}  \ar[lu]
\ar[ld]
\\ {} & {\bullet} & {\diamond} \ar[lu]\ar[ld] & {\rhd} \ar[lu] \ar[ld] \\ {} & {\bullet}
& {\diamond} \ar[lu] \ar[ld] & {\rhd} \ar[lu] \ar[ld] \\
{} & {\bullet} & {\diamond} \ar[lu] \ar[ld] & {\rhd} \ar[lu]
\ar[ld] \\ {} & ... & ...\ar[lu] & ...\ar[lu]} \ee where every
$\diamond$ (resp. $\rhd$) is obtained from $\bullet$ (resp.
$\diamond$) of the same generalized weight by applying the
transpose of $h(0)$. For simplicity we do not display doted arrows
and arrows obtained by the "addition of arrows" rule.
\end{corollary}
\noindent {\em Proof.} From an embedding $M(1) \otimes \Omega_2
\hookrightarrow M(1) \otimes \Omega_3$ it is clear that there will
be arrows connecting $u^{2,m}$ and $u^{m \pm 1}$ as in Theorem
\ref{2step}. Let $u^{3,n}=h^T(0)u^{2,n}$, so $u^{3,n}$ is
represented by a $\rhd$. By arguing as before, from each $u^{3,m}$
there will be arrows pointing  to $u^{2,m+1}$ and $u^{2,m-1}$. So
we only have to show that there are no additional arrows from
$u^{3,m}$ except those displayed on (\ref{embed2}). From the
formula $L(0)u^{3,n}=n^2 u^{3,n}+\frac{1}{2}u^n$ it follows that
$u^n$ and $u^{3,n}$ form a Jordan block with respect to $L(0)$.
But $u^{n}$ can be also reached from $u^{3,n}$ via an oriented
path $u^{3,n} \rightarrow u^{2,n+1} \rightarrow u^n$ so there is
no need to display an arrow from $u^{3,n}$ to $u^n$, because the
submodule generated by $u^{3,n}$ contains $u^{n}$. Similarly with
$u^{3,n}$ and $u^{n+2}$.
%Recall that $L(n)$, $n \neq 0$ is an infinite sum of products of
%two Heisenberg algebra generators, so it is impossible to "jump"
%with $L(n)$, $n \neq 0$ from $\rhd$ to $\bullet$ without "passing
%through" an intermediate $\diamond$.
%Notice that the
%identity $[L(1),L(-1)]=2L(0)$ represents "addition of vectors"
%$\nwarrow+\swarrow=\longleftarrow$ in (\ref{embed2}).
\epfv

%\section{"Visualizing" logarithmic modules: Systems with regular singularities }
%We assume some knowledge of vertex operator algebra theory as in
%\cite{FLM}, \cite{K} and \cite{FrB}). Also we assume that the reader
%is familiar with the basics of the theory of intertwining operators as in
%the last chapter of \cite{FHL}.
%{\em Facts of Life:} \\
%Suppose that we have a first-order system with regular
%singularities
%$$\frac{dY}{dx}=AY,$$
%where $A \mapsto \mathbb{C}^* \longrightarrow M_n(\mathbb{C})$ is
%a holomorphic mapping on the punctured disk. The it is known that
%
%$$\langle w'_4, \mathcal{Y}(w_1,z_1)\mathcal{Y}(w_2,z_2)w_3 \rangle.$$
%It is known that we can write a differential equation with regular
%singularities.

\section{Hidden logarithmic intertwining operators}

%In the previous section we gave a complete description of $M(1)
%\otimes \Omega$ and $M(1) \otimes \Omega_3$in the case when
%$h(0)|_{\Omega_2}$ is nilpotent. By using Theorem \ref{griess} we
%first construct $\mathcal{Y} \in \ I \ { L(1,i^2) \choose L(1,m^2)
%\ L(1,n^2)}$, where $i \in \{ m+n,m+n-2,...,|m-n| \}$. Then we can
%construct intertwining operators of the type
%For instance,consider two indecomposable module
%$$W_3 \choose  \xymatrix{{\bullet} &  & \\ & \ar[ul] {\diamond}} \ \ \ \xymatrix{{\bullet} & \\ & \ar[ul] {\diamond}}
%$$
%%$$\xymatrix{ & \ar[ld] {\diamond} \\ {\bullet} & }$$
%Let $\Omega$ be a two-dimensional space.
Suppose that $h(0)|_{\Omega_2}$ in some basis $\{w_1,w_2\}$ for
$\Omega_2$ is represented by $\left[\begin{array}{cc} 0 & 1
\\ 0 & 0
\end{array}
\right]$, and $h(0)|_{\Omega_3}$, in some basis
$\{\tilde{w}_1,\tilde{w}_2,\tilde{w}_3\}$ for $\Omega_3$,
is represented by $\left[ \begin{array}{ccc} 0 & 1 & 0 \\
0 & 0 & 1
\\ 0 & 0 & 0 \end{array} \right]$.
Then a surjective map $$T_{\Omega_2,\Omega_2}^{\Omega_3} \in {\rm
Hom}_{}( \Omega_2 \otimes \Omega_2, \Omega_3),$$ defined by \be
\label{mapt} w_2 \otimes w_2 \mapsto \tilde{w}_3, \ w_2 \otimes
w_1 \mapsto \frac{\tilde{w}_2}{2}, \ {w_1 \otimes w_2} \mapsto
\frac{\tilde{w}_2}{2} , \ w_1 \otimes w_1 \mapsto
\frac{\tilde{w}_1}{2} \ee commutes with $h(0)$.

Let us denote by $W_2(1,m^2) \subset M(1) \otimes \Omega_2$ a
cyclic $Vir$-module generated by $u^{2,m}$ of weight $m^2$. For $m
> 0$, $W_2(1,m^2)$ can be visualized as a "wedge" in (\ref{embed})
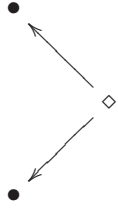
$$\xymatrix{{\bullet} &  \\
& \diamond \ar[ul] \ar[ld] \\ \bullet &}$$ or a single arrow
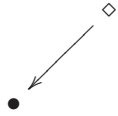
$$\xymatrix{& \diamond \ar[ld] \\
\bullet & }$$ in the $m=0$ case. Similarly, we denote by
$W_3(1,m^2) \subset M(1) \otimes \Omega_3$ (cf. Corollary
\ref{om-3}) the module generated by $\rhd$, of generalized weight
$m^2$. For every $m >1$ this module can be visualized as
$$
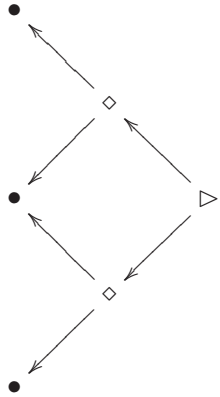
\xymatrix{{\bullet} & &  \\
&  \diamond \ar[lu] \ar[ld] & \\ \bullet & & \rhd \ar[lu] \ar[ld] \\ & \diamond \ar[ul] \ar[ld] &  \\
\bullet & &}$$ Similarly, $W_3(1,0)$ may be visualized as
$$
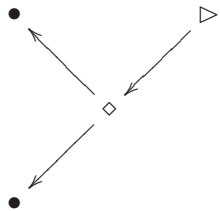
\xymatrix{ {\bullet} & & \rhd \ar[ld]  \\
& \diamond \ar[lu] \ar[ld] &  \\ \bullet & }$$ Again, we shall
assume that $\rhd$ (resp. $\diamond$) is obtained from $\diamond$
(resp. $\bullet$) of the same generalized weight by applying
$h(0)^T$. The following Lemma is just a consequence of
$L(0)u^{3,m}=m^2 u^{3,m}+\frac{1}{2}u^{m}$, so we omit the proof.
\begin{lemma} \label{genuine}
For every $m \geq 0$ the module $W_3(1,m^2)$ is a genuine
logarithmic module.
\end{lemma}
The module $W_3(1,m^2)$ is a nonsplit extension of $L(1,m^2)$ by
$W_2(1,(m-1)^2) + W_2(1,(m+1)^2)$, where if $m=0$ the first
summand is trivial.
%The irreducible submodules of $W_3(1,m^2)$ can
%be easily classified as in the case of $W_2(1,m^2)$. In
%particular, every proper submodule of $W_3(1,m^2)$ is contained
%inside $W_2(1,(m-1)^2) +W_2(1,(m+1)^2)$.

\noindent Now we have a logarithmic version of Theorem
\ref{griess}.
\begin{theorem} \label{griesslog}
Let $\Omega_2$ and $\Omega_3$ be as above. Then there exists a
nontrivial $\mathcal{Y} \in I \ \ds{{ M(1) \otimes \Omega_3
\choose M(1) \otimes \Omega_2 \ M(1) \otimes \Omega_2 }}$, such
that
%$T_{\Omega_1,\Omega_2}^{\Omega_3}$.
$\mathcal{Y}$ projects down to a hidden logarithmic intertwining
operator
$$\bar{\mathcal{Y}} \in I \ {W \choose W_2(1,m^2) \ W_2(1,n^2) }$$
of depth one, where
\begin{equation} \label{wdirect}
W=\sum_{{\tiny \begin{array}{c} |m-n| \leq k \leq m+n \\ k \equiv
m+n \ mod \ 2 \end{array}}} W_3(1,k^2).
\end{equation}
This sum is not direct, whenever $mn \neq 0$.
\end{theorem}
\noindent {\em Proof.} Let $\mathcal{Y}$ be as in Theorem
\ref{pre-upper}, with $T_{\Omega_2,\Omega_2}^{\Omega_3}$ as in
(\ref{mapt}). Let
$$\bar{\mathcal{Y}}(\cdot \ ,x) \cdot \ =\mathcal{Y}( \cdot \ ,x) \cdot \ |_{W_2(1,m^2)
\otimes W_2(1,n^2)}.$$ We recall (cf. Section 7) that in this case
\bea &&
\mathcal{Y}(u^0,x)=\mathcal{Y}(w_1,x)=T_{\Omega_2,\Omega_2}^{\Omega_3}(w_1),
\nonumber \\
\label{mnk} && \mathcal{Y}(u^{2,0},x)=\mathcal{Y}(w_2,x)= \int^-
h(x)
T_{\Omega_2,\Omega_2}^{\Omega_3}(w_1)+T_{\Omega_2,\Omega_2}^{\Omega_3}(w_1)\int^+
h(x) +T_{\Omega_2,\Omega_2}^{\Omega_3}(w_2). \eea Also, \be
\label{mnk1} \mathcal{Y}(u^{2,m},x)=\mathcal{Y}(P_m(h) u^{2,0},x).
\ee We need a more precise information about the image of
$\mathcal{Y}$. Since $W_2(1,m^2)$ is cyclic, and generated by
$u^{2,m}$, the image of $\bar{\mathcal{Y}}$, denoted by $\hat{W}$,
is actually the Virasoro submodule generated by the Fourier
coefficients of $\mathcal{\bar{Y}}(u^{2,m},x)u^{2,n}$,
$\mathcal{\bar{Y}}(u^{m \pm 1},x)u^{2,n}$,
$\mathcal{\bar{Y}}(u^{2,m},x)u^{n \pm 1}$ and
$\mathcal{\bar{Y}}(u^{m \pm 1},x)u^{n \pm 1}$. As before, let
$u^{3,k}$ denote a generator of $W_3(1,k^2)$ of generalized weight
$k^2$.
%Then $\mathcal{Y}$ is of type $\ds{{M(1)\otimes \Omega_3 \choose
%M(1) \otimes \Omega_2 \ M(1) \otimes \Omega_2}}$.
Since $L(1,(m-1)^2) \oplus L(1,(m+1)^2) \subset W_2(1,m^2)$, then
Theorem \ref{griess}, (\ref{mnk}) and (\ref{mnk1}) imply that the
submodule of $\hat{W}$ generated by $\bar{\mathcal{Y}}(u^{m \pm
1},x)u^{n \pm 1}$ is precisely \be \label{ssss} L(1,(m-n-2)^2)
\oplus L(1,(m-n)^2) \oplus \cdots \oplus L(1,(m+n+2)^2). \ee Now
we move "one step higher" or " deeper" in the filtration and
determine the Virasoro submodule generated by
$\bar{\mathcal{Y}}(u^{2,m},x)u^n$. From the formula
$$h(0) {\mathcal{Y}} (u^{2,m},x)u^n=\mathcal{Y}(u^m,x)u^n,$$
and the previous discussion it follows that the submodule
generated by $\mathcal{Y}(u^{2,m},x)u^{n \pm 1}$ is contained
inside $W_2(1,(m-n-1)^2) + \cdots + W_2(1,(m+n+1)^2)$ and possibly
some $L(1,k^2)$, $k^2 \notin \{ (m-n-2)^2,...,(m+n+2)^2\}$. But
having in the image such $L(1,k^2)$ would contradict to Theorem
\ref{griess}. Similarly, the submodule generated by
$\bar{\mathcal{Y}}(u^{m \pm 1},x)u^{2,n}$ lies again inside the
sum $W_2(1,(m-n-1)^2) + \cdots + W_2(1,(m+n-1)^2)$. Furthermore,
from
$$h(0) {\mathcal{Y}}(u^{2,m},x)u^{2,n}={\mathcal{Y}}(u^m,x)u^{2,n}+{\mathcal{Y}}(u^{2,m},x)u^n$$
it follows that the Virasoro module generated by the coefficients
of $\bar{\mathcal{Y}}(u^{2,m},x)u^{2,n}$ is contained inside
$W_3(1,(m-n)^2) + \cdots + W_3(1,(m+n)^2)$ and possibly some
irreducible module $L(1,k^2)$ not included in (\ref{ssss}). But
this would again contradict to Theorem \ref{griess}. Thus, we have
shown that the image $\hat{W}$ is contained inside $W$ (cf.
\ref{wdirect}). In fact, it is not hard to show that $u^{3,k} \in
\hat{W}$ for $k^2 \in \{(m-n)^2,...,(m+n)^2\}$, which would imply
$\hat{W}=W$.
% In
%fact, from the discussion following Lemma \ref{genuine} it
%suffices to show that
%$$u^{3,k}+u_k \in W, \ \ u_k \in W_2(1,(m-n-1)^2) + \cdots +
%W_2(1,(m+n+1)^2),$$ for every $k$ as above. But this follows from
%the definition of (\ref{mnk}) and (\ref{mnk1}).
%(cf. Section 7) that \bea \label{mnk} &&
%\mathcal{Y}(u^{2,0},x)=\mathcal{Y}(w_2,x)= \nn && \int^- h(x)
%T_{\Omega_2,\Omega_2}^{\Omega_3}(w_1)+T_{\Omega_2,\Omega_2}^{\Omega_3}(w_1)\int^+
%h(x) +T_{\Omega_2,\Omega_2}^{\Omega_3}(w_2). \eea The previous
%formulas combined with
%$$\mathcal{Y}(u^{2,m},x)=\mathcal{Y}(P_m(h) u^{2,0},x),$$
%and the definition of $T_{\Omega_2,\Omega_2}^{\Omega_3}$ proves
%the claim.

Finally, from (\ref{mnk}) and (\ref{mnk1}) it is clear that
$\bar{\mathcal{Y}}$ is a genuine logarithmic intertwining operator
of depth one. \epfv

%More generally,
%$$\mathcal{Y}(h(-i_1-1) \cdots h(-i_k-1) \cdot u^{2,0},x)= \nb \frac{d^{i_1}}{dx} h(x) \cdots \frac{d^{i_1}}{dx} h(x)
%T_{\Omega_2,\Omega_2}^{\Omega_3}(w_1) \int h(x) \nb $$ $$+ \nb
%\frac{d^{i_1}}{dx} h(x) \cdots \frac{d^{i_k}}{dx} h(x) \nb
%T_{\Omega_2,\Omega_2}^{\Omega_3}(w_2).$$
The module $W$ in the previous theorem has the following diagram
representation
$$
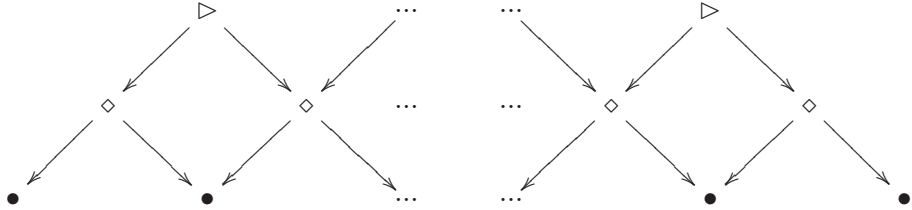
\xymatrix{ & & & \rhd \ar[dr] \ar[dl] & & ...\ar[ld] & ...\ar[rd]
& & \rhd \ar[dl] \ar[dr] & &
\\ & & \diamond \ar[dr] \ar[dl] &&
\diamond \ar[dr]\ar[dl] & ... & ... & \diamond \ar[rd]\ar[ld] &&
\diamond \ar[dr] \ar[dl]
& \\
& \bullet & & \bullet && ... & ... &  & \bullet && \bullet }
$$
%the corresponding intertwining operator can be symbolically
%visualized as being of type {\tiny
%$$\left( \begin{array}{c}
%\xymatrix{{\bullet} & &  \\
%&  \diamond \ar[lu] \ar[ld] & \\ \bullet & & \rhd \ar[lu] \ar[ld] \\ & \diamond \ar[ul] \ar[ld] &  \\
%\bullet & &} \\
%\xymatrix{{\bullet} & & \\
%& \diamond \ar[ul] \ar[ld] & \\ \bullet & & } \ \ \
%\xymatrix{{\bullet} & & \\
%& \diamond \ar[ul] \ar[ld] &  \\ \bullet & &}
%\end{array}
%\right)
%$$
%}

\begin{remark}
{\em It is not hard to generalize the results from this section to
logarithmic modules with Jordan blocks of arbitrary size.
%Another direction is to find explicitly description of the fusion
%algebra among $W_n(1,m^2)$ and their submodules by combining ideas
%from this paper \cite{M1} and \cite{M2}
}
\end{remark}

\section{Hidden logarithmic intertwining operators among Feigin-Fuchs modules at $c=1-12a^2$}

Let us recall that $M(1,\lambda)_a$ is a $M(1)_a$-module of lowest
conformal weight $\frac{\lambda^2}{2}-a \lambda$. As we have
already mentioned in the introduction, the $\lambda=a$ case (e.g.,
$\lambda=0$ for $a=0$) is indeed very special. Here is a
consequence of Corollary \ref{hha}.
\begin{lemma} \label{333}
Let $M(1)_a \otimes \Omega$, where $\Omega$ is two-dimensional and
$h(0)|_{\Omega}$ is represented by $\left[
\begin{array}{cc} \lambda & 1 \\ 0 & \lambda
\end{array}\right]$, in some basis. Then $M(1)_a \otimes \Omega$ is an ordinary $M(1)_a$-module if and only if
$\lambda=a$, in which case the Feigin-Fuchs module $M(1,a)_a$ is
of lowest conformal weight $\frac{-a^2}{2}$.
\end{lemma}
Let $\Omega$ be as in the lemma. It is easy to see that $M(1)_a
\otimes (\Omega \otimes \Omega)$ is a genuine logarithmic
$M(1)_a$-modules. Now, we have a consequence of Theorem
\ref{pre-upper}.
\begin{corollary} \label{last-cor}
Let $M(1)_a \otimes \Omega$ be as in Lemma \ref{333}. Then there
exists a genuine logarithmic intertwining operator of type $\ds{{
M(1)_a \otimes (\Omega \otimes \Omega) \choose M(1)_a \otimes
\Omega \ \ M(1)_a \otimes \Omega}}$.
\end{corollary}

\noindent {\bf Example.} For $a=\frac{1}{2}$ the vertex operator
algebra $M(1)_{\frac{1}{2}}$ has central charge $c=-2$. Then
$M(1)_{\frac{1}{2}}$-module $M(1)_{\frac{1}{2}} \otimes \Omega$,
where $h(0)|_{\Omega}$ is represented by
$$\left[ \begin{array}{cc} \frac{1}{2} & 1 \\ 0 & \frac{1}{2} \end{array}
\right]$$ is an ordinary $M(1)_{\frac{1}{2}}$-module with the
lowest conformal weight $-\frac{1}{8}$. The intertwining operator
constructed in Corollary \ref{last-cor} is closely related to
logarithmic operators studied in \cite{Gu}.

% Then,
%$$L(0)|_{\Omega}=(\frac{h(0)^2-h(0)}{2})_{\Omega}=\left[\begin{array}{cc} \frac{-1}{8} & 0 \\ 0 & \frac{-1}{8} \end{array}
%\right] $$

%\section{Future work}

%In the sequel we will study in more details various subalgebras
%of $M(1)_a$ and associated logarithmic intertwining operators
%\cite{M3}. Many results from this paper can be extended to other
%higher rank $\mathcal{W}$-algebras by using the construction in
%\cite{FKRW}. In \cite{AM} we shall describe hidden intertwining
%operators for $\mathcal{W}(2,2p-1)$-algebras with two generators.

\end{document}